\documentclass{amsproc}
\usepackage{amsmath}
\usepackage{amssymb,latexsym}
\usepackage{amscd}
\usepackage{xspace}
\usepackage{verbatim}
\begin{document}
\numberwithin{equation}{section}
\title[Precise boundary trace]{The precise boundary trace of positive solutions of the equation $\Gd u=u^q$ in the supercritical case.}
\author{Moshe Marcus }
\address{Department of Mathematics, Technion\\
 Haifa 32000, ISRAEL}
 \email{marcusm@math.technion.ac.il}
\author{ Laurent Veron}
\address{Laboratoire de Math\'ematiques, Facult\'e des Sciences\\
Parc de Grandmont, 37200 Tours, FRANCE}
\email{veronl@lmpt.univ-tours.fr}
\date{\today}

\newcommand{\txt}[1]{\;\text{ #1 }\;}
\newcommand{\tbf}{\textbf}
\newcommand{\tit}{\textit}
\newcommand{\tsc}{\textsc}
\newcommand{\trm}{\textrm}
\newcommand{\mbf}{\mathbf}
\newcommand{\mrm}{\mathrm}
\newcommand{\bsym}{\boldsymbol}
\newcommand{\scs}{\scriptstyle}
\newcommand{\sss}{\scriptscriptstyle}
\newcommand{\txts}{\textstyle}
\newcommand{\dsps}{\displaystyle}
\newcommand{\fnz}{\footnotesize}
\newcommand{\scz}{\scriptsize}
\newcommand{\be}{\begin{equation}}
\newcommand{\bel}[1]{\begin{equation}\label{#1}}
\newcommand{\ee}{\end{equation}}
\newtheorem{subn}{\name}
\renewcommand{\thesubn}{}
\newcommand{\bsn}[1]{\def\name{#1$\!\!$}\begin{subn}}
\newcommand{\esn}{\end{subn}}
\newtheorem{sub}{\name}[section]
\newcommand{\dn}[1]{\def\name{#1}}   
\newcommand{\bs}{\begin{sub}}
\newcommand{\es}{\end{sub}}
\newcommand{\bsl}[1]{\begin{sub}\label{#1}}
\newcommand{\bth}[1]{\def\name{Theorem}\begin{sub}\label{t:#1}}
\newcommand{\blemma}[1]{\def\name{Lemma}\begin{sub}\label{l:#1}}
\newcommand{\bcor}[1]{\def\name{Corollary}\begin{sub}\label{c:#1}}
\newcommand{\bdef}[1]{\def\name{Definition}\begin{sub}\label{d:#1}}
\newcommand{\bprop}[1]{\def\name{Proposition}\begin{sub}\label{p:#1}}
\newcommand{\bnote}[1]{\def\name{\mdseries\scshape Notation}\begin{sub}\label{n:#1}}
\newcommand{\bproof}{\begin{proof}}
\newcommand{\eproof}{\end{proof}}
\newcommand{\bcom}{}
\newcommand{\req}{\eqref}
\newcommand{\rth}[1]{Theorem~\ref{t:#1}}
\newcommand{\rlemma}[1]{Lemma~\ref{l:#1}}
\newcommand{\rcor}[1]{Corollary~\ref{c:#1}}
\newcommand{\rdef}[1]{Definition~\ref{d:#1}}
\newcommand{\rprop}[1]{Proposition~\ref{p:#1}}
\newcommand{\rnote}[1]{Notation~\ref{n:#1}}
\newcommand{\BA}{\begin{array}}
\newcommand{\EA}{\end{array}}
\newcommand{\BAN}{\renewcommand{\arraystretch}{1.2}
\setlength{\arraycolsep}{2pt}\begin{array}}
\newcommand{\BAV}[2]{\renewcommand{\arraystretch}{#1}
\setlength{\arraycolsep}{#2}\begin{array}}
\newcommand{\BSA}{\begin{subarray}}
\newcommand{\ESA}{\end{subarray}}
\newcommand{\BAL}{\begin{aligned}}
\newcommand{\EAL}{\end{aligned}}
\newcommand{\BALG}{\begin{alignat}}
\newcommand{\EALG}{\end{alignat}}
\newcommand{\BALGN}{\begin{alignat*}}
\newcommand{\EALGN}{\end{alignat*}}
\newcommand{\note}[1]{\noindent\textit{#1.}\hspace{2mm}}
\newcommand{\Remark}{\note{Remark}}
\newcommand{\forevery}{\quad \forall}
\newcommand{\1}{\\[1mm]}
\newcommand{\2}{\\[2mm]}
\newcommand{\3}{\\[3mm]}
\newcommand{\set}[1]{\{#1\}}
\def\({{\rm (}}
\def\){{\rm )}}
\newcommand{\st}[1]{{\rm (#1)}}
\newcommand{\lra}{\longrightarrow}
\newcommand{\lla}{\longleftarrow}
\newcommand{\llra}{\longleftrightarrow}
\newcommand{\Lra}{\Longrightarrow}
\newcommand{\Lla}{\Longleftarrow}
\newcommand{\Llra}{\Longleftrightarrow}
\newcommand{\warrow}{\rightharpoonup}
\def\dar{\downarrow}
\def\uar{\uparrow}
\newcommand{\paran}[1]{\left (#1 \right )}
\newcommand{\sqrbr}[1]{\left [#1 \right ]}
\newcommand{\curlybr}[1]{\left \{#1 \right \}}
\newcommand{\abs}[1]{\left |#1\right |}
\newcommand{\norm}[1]{\left \|#1\right \|}
\newcommand{\angbr}[1]{\left< #1\right>}
\newcommand{\paranb}[1]{\big (#1 \big )}
\newcommand{\sqrbrb}[1]{\big [#1 \big ]}
\newcommand{\curlybrb}[1]{\big \{#1 \big \}}
\newcommand{\absb}[1]{\big |#1\big |}
\newcommand{\normb}[1]{\big \|#1\big \|}
\newcommand{\angbrb}[1]{\big\langle #1 \big \rangle}
\newcommand{\thkl}{\rule[-.5mm]{.3mm}{3mm}}
\newcommand{\thknorm}[1]{\thkl #1 \thkl\,}
\newcommand{\trinorm}[1]{|\!|\!| #1 |\!|\!|\,}
\newcommand{\vstrut}[1]{\rule{0mm}{#1}}
\newcommand{\rec}[1]{\frac{1}{#1}}
\newcommand{\opname}[1]{\mathrm{#1}\,}
\newcommand{\supp}{\opname{supp}}
\newcommand{\dist}{\opname{dist}}
\newcommand{\sign}{\opname{sign}}
\newcommand{\diam}{\opname{diam}}
\newcommand{\q}{\quad}
\newcommand{\qq}{\qquad}
\newcommand{\hsp}[1]{\hspace{#1mm}}
\newcommand{\vsp}[1]{\vspace{#1mm}}
\newcommand{\prt}{\partial}
\newcommand{\sms}{\setminus}
\newcommand{\ems}{\emptyset}
\newcommand{\ti}{\times}
\newcommand{\pr}{^\prime}
\newcommand{\ppr}{^{\prime\prime}}
\newcommand{\tl}{\tilde}
\newcommand{\wtl}{\widetilde}
\newcommand{\sbs}{\subset}
\newcommand{\sbeq}{\subseteq}
\newcommand{\nind}{\noindent}
\newcommand{\ovl}{\overline}
\newcommand{\unl}{\underline}
\newcommand{\nin}{\not\in}
\newcommand{\pfrac}[2]{\genfrac{(}{)}{}{}{#1}{#2}}
\newcommand{\tin}{\to\infty}
\newcommand{\ind}[1]{_{_{#1}}\!}
\newcommand{\chr}[1]{\chi\ind{#1}}
\newcommand{\rest}[1]{\big |\ind{#1}}
\newcommand{\num}[1]{{\rm (#1)}\hspace{2mm}}
\newcommand{\wkc}{weak convergence\xspace}
\newcommand{\wrto}{with respect to\xspace}
\newcommand{\cons}{consequence\xspace}
\newcommand{\consy}{consequently\xspace}
\newcommand{\Consy}{Consequently\xspace}
\newcommand{\Essy}{Essentially\xspace}
\newcommand{\essy}{essentially\xspace}
\newcommand{\mnz}{minimizer\xspace}
\newcommand{\sth}{such that\xspace}
\newcommand{\ngh}{neighborhood\xspace}
\newcommand{\nghs}{neighborhoods\xspace}
\newcommand{\seq}{sequence\xspace}
\newcommand{\sseq}{subsequence\xspace}
\newcommand{\ifif}{if and only if\xspace}
\newcommand{\suff}{sufficiently\xspace}
\newcommand{\abc}{absolutely continuous\xspace}
\newcommand{\sol}{solution\xspace}
\newcommand{\subsol}{subsolution\xspace}
\newcommand{\supsol}{supersolution\xspace}
\newcommand{\Wlg}{Without loss of generality\xspace}
\newcommand{\wlg}{without loss of generality\xspace}
\newcommand{\bdw}{\partial\Gw}
\newcommand{\Capq}{C_{2/q,q'}}
\newcommand{\Wq}{W^{2/q,q'}}
\newcommand{\Wqdual}{W^{-2/q,q}}
\newcommand{\Wqdb}{W^{-2/q,q}_+(\bdw)}
\newcommand{\sbsq}{\overset{q}{\sbs}}
\newcommand{\smq}{\overset{q}{\sim}}
\newcommand{\app}[1]{\underset{#1}{\approx}}
\newcommand{\suppq}{\mathrm{supp}^q_{\bdw}\,}
\newcommand{\convq}{\overset{q}{\to}}
\newcommand{\barq}[1]{\bar{#1}^{^q}}
\newcommand{\prtq}{\partial_q}
\newcommand{\tr}{\mathrm{tr}\,}
\newcommand{\Tr}{\mathrm{Tr}\,}
\newcommand{\trR}{\mathrm{tr}\ind{\CR}}
\newcommand{\trin}[1]{\mathrm{tr}\ind{#1}}
\newcommand{\qcl}{$q$-closed\xspace}
\newcommand{\qop}{$q$-open\xspace}
\newcommand{\gsmod}{$\gs$-moderate\xspace}
\newcommand{\gsreg}{$\gs$-regular\xspace}
\newcommand{\qreg}{$q$-quasi regular\xspace}
\newcommand{\qeq}{$q$-equivalent\xspace}
\newcommand{\ppf}{\underset{f}{\prec\prec}}
\newcommand{\ofrown}{\overset{\frown}}
\newcommand{\modcon}{\underset{mod}{\lra}}
\newcommand{\ugb}[1]{u\chr{\Gs_\gb(#1)}}
\newcommand{\mcon}{$q$-moderately convergent\xspace}
\newcommand{\mdiv}{$q$-moderately divergent\xspace}
\def\qsupp{q\text{-supp}\,}
\def\Lim{\,\text{\rm Lim}\,}
\def\muCR{\mu\ind{\CR}}
\def\vCR{v\ind{\CR}}
\def\bcom{}
\def\ga{\alpha}     \def\gb{\beta}       \def\gg{\gamma}
\def\gc{\chi}       \def\gd{\delta}      \def\ge{\epsilon}
\def\gth{\theta}                         \def\vge{\varepsilon}
\def\gf{\phi}       \def\vgf{\varphi}    \def\gh{\eta}
\def\gi{\iota}      \def\gk{\kappa}      \def\gl{\lambda}
\def\gm{\mu}        \def\gn{\nu}         \def\gp{\pi}
\def\vgp{\varpi}    \def\gr{\rho}        \def\vgr{\varrho}
\def\gs{\sigma}     \def\vgs{\varsigma}  \def\gt{\tau}
\def\gu{\upsilon}   \def\gv{\vartheta}   \def\gw{\omega}
\def\gx{\xi}        \def\gy{\psi}        \def\gz{\zeta}
\def\Gg{\Gamma}     \def\Gd{\Delta}      \def\Gf{\Phi}
\def\Gth{\Theta}
\def\Gl{\Lambda}    \def\Gs{\Sigma}      \def\Gp{\Pi}
\def\Gw{\Omega}     \def\Gx{\Xi}         \def\Gy{\Psi}

\def\CS{{\mathcal S}}   \def\CM{{\mathcal M}}   \def\CN{{\mathcal N}}
\def\CR{{\mathcal R}}   \def\CO{{\mathcal O}}   \def\CP{{\mathcal P}}
\def\CA{{\mathcal A}}   \def\CB{{\mathcal B}}   \def\CC{{\mathcal C}}
\def\CD{{\mathcal D}}   \def\CE{{\mathcal E}}   \def\CF{{\mathcal F}}
\def\CG{{\mathcal G}}   \def\CH{{\mathcal H}}   \def\CI{{\mathcal I}}
\def\CJ{{\mathcal J}}   \def\CK{{\mathcal K}}   \def\CL{{\mathcal L}}
\def\CT{{\mathcal T}}   \def\CU{{\mathcal U}}   \def\CV{{\mathcal V}}
\def\CZ{{\mathcal Z}}   \def\CX{{\mathcal X}}   \def\CY{{\mathcal Y}}
\def\CW{{\mathcal W}}
\def\BBA {\mathbb A}   \def\BBb {\mathbb B}    \def\BBC {\mathbb C}
\def\BBD {\mathbb D}   \def\BBE {\mathbb E}    \def\BBF {\mathbb F}
\def\BBG {\mathbb G}   \def\BBH {\mathbb H}    \def\BBI {\mathbb I}
\def\BBJ {\mathbb J}   \def\BBK {\mathbb K}    \def\BBL {\mathbb L}
\def\BBM {\mathbb M}   \def\BBN {\mathbb N}    \def\BBO {\mathbb O}
\def\BBP {\mathbb P}   \def\BBR {\mathbb R}    \def\BBS {\mathbb S}
\def\BBT {\mathbb T}   \def\BBU {\mathbb U}    \def\BBV {\mathbb V}
\def\BBW {\mathbb W}   \def\BBX {\mathbb X}    \def\BBY {\mathbb Y}
\def\BBZ {\mathbb Z}

\def\GTA {\mathfrak A}   \def\GTB {\mathfrak B}    \def\GTC {\mathfrak C}
\def\GTD {\mathfrak D}   \def\GTE {\mathfrak E}    \def\GTF {\mathfrak F}
\def\GTG {\mathfrak G}   \def\GTH {\mathfrak H}    \def\GTI {\mathfrak I}
\def\GTJ {\mathfrak J}   \def\GTK {\mathfrak K}    \def\GTL {\mathfrak L}
\def\GTM {\mathfrak M}   \def\GTN {\mathfrak N}    \def\GTO {\mathfrak O}
\def\GTP {\mathfrak P}   \def\GTR {\mathfrak R}    \def\GTS {\mathfrak S}
\def\GTT {\mathfrak T}   \def\GTU {\mathfrak U}    \def\GTV {\mathfrak V}
\def\GTW {\mathfrak W}   \def\GTX {\mathfrak X}    \def\GTY {\mathfrak Y}
\def\GTZ {\mathfrak Z}   \def\GTQ {\mathfrak Q}
\font\Sym= msam10
\def\SYM#1{\hbox{\Sym #1}}
\def\rqq{\req{q-eq}\xspace}

\subjclass{Primary 35J60, 35J67; Secondary 31B15, 31B20, 31C15}

\dedicatory{To Ha\"\i m, with friendship and high esteem.}

\keywords{Nonlinear elliptic equations, Bessel capacities, fine
topology, Borel measures}

\begin{abstract}
We construct the precise boundary trace of positive solutions of
$\Delta u=u^q$ in a smooth bounded domain $\Gw\sbs\BBR^N$, for $q$
in the super-critical case $q\geq (N+1)/(N-1)$. The construction
is performed in the framework of the fine topology associated with
the Bessel capacity $C_{{2/q,q'}}$ on $\bdw$. We prove that the
boundary trace is a Borel measure (in general unbounded),which is
outer regular and essentially absolutely continuous relative to
this capacity. We provide a necessary and sufficient condition for
such measures to be the boundary trace of a positive solution and
prove that the corresponding generalized boundary value problem is
well-posed in the class of $\sigma$-moderate solutions.
\end{abstract}

\maketitle

\section{Introduction} In this paper we present a theory of boundary trace of positive solutions of the equation
\begin{equation}\label{q-eq}
 -\Gd u+|u|^{q-1}u=0
\end{equation}
 in a  bounded domain $\Gw\sbs \BBR^N$ of class $C^2$. A function $u$ is a solution if
 $u\in L^q_{loc}(\Gw)$ and the equation holds in the distribution sense.
\par Semilinear elliptic equations with absorption, of which \rqq is one of the most important,
have been intensively studied in the last 30 years. The foundation
for these studies can be found in the pioneering work of Brezis
starting with his joint research with Benilan in the 70's
\cite{BB}, and followed by a series of works with colleagues and
students, up to the present.
\par In the subcritical case, $1<q<q_c=(N+1)/(N-1)$, the boundary trace
 theory and the associated boundary value problem, are well understood. This theory has
  been developed, in parallel, by two different methods: one based on a probabilistic approach
 (see \cite{LG, DK1, DK2}, Dynkin's book \cite{Dbook1} and the references therein) and the other
purely analytic (see \cite{MV1,MV2,MV3}). In 1997 Le Gall showed that this theory is not appropriate for the supercritical
case because, in this case, there may be infinitely many solutions with the same boundary trace.
Following this observation, a  theory of 'fine' trace was introduced by
Dynkin and Kuznetsov  \cite{DK3}. Their results demonstrated that, for $q\leq 2$,
the fine trace theory
 is satisfactory in the family of so-called \gsmod solutions. A few years later Mselati \cite{Ms}
 used this theory and other results of Dynkin \cite{Dbook1}, in combination with the Brownian snake method
 developed by Le Gall \cite{LGbook}, in order to show that, in the case $q=2$ all positive solutions
 are \gsmod. Shortly thereafter Marcus and Veron
\cite{MV4} proved that, for all $q\geq q_c$ and every compact set
$K\sbs \bdw$,
 the maximal solution of \rqq vanishing outside $K$ is \gsmod. Their proof
was based on the derivation of sharp capacitary estimates for the
maximal solution.
 In continuation, Dynkin \cite{Dbook2} used
Mselati's (probabilistic) approach and the results of  Marcus and
Veron \cite{MV4}  to show that, in the case $q\leq 2$, all
positive solutions are \gsmod.
 For $q>2$ the problem remains open.

Our definition of boundary trace is based on the fine topology
associated
 with the Bessel capacity $\Capq$ on $\bdw$, denoted by $\GTT_q$.
The presentation requires some notation. \1
\note{Notation 1.1}
\begin{description}
  \item[a]   For every $x\in \BBR^N$ and every $\gb>0$ put
  $\gr(x):=\dist(x,\prt\Gw)$ and
$$\Gw_\gb=\set{x\in \Gw:\,\gr(x)<\gb},\q \Gw'_\gb=\Gw\sms \bar \Gw_\gb, \q \Gs_\gb=\prt\Gw'_\gb.$$
  \item[b]  There exists a positive number $\gb_0$ \sth,
\begin{equation}\label{gsgr}
  \forall x\in \Gw_{\gb_0}\q \exists !\; \gs(x)\in\bdw \,:\, \dist(x,\gs(x))=\gr(x).
\end{equation}
If (as we assume) $\Gw$ is of class $C^2$ and
$\gb_0$ is sufficiently small, the mapping $x\mapsto (\gr(x),\gs(x))$ is a $C^2$ diffeomorphism of $\Gw_{\gb_0}$
onto $(0,\gb_0)\ti\bdw$.
\item[c] If $Q\sbs \bdw$  put $\Gs_\gb(Q)=\{x\in \Gs_\gb:\, \gs(x)\in Q\}$.
\item[d] If $Q$ is a $\GTT_q$-open subset of $\bdw$ and $u\in C(\bdw)$ we denote by $u_\gb^Q$
the solution of \rqq in $\Gw'_\gb$ with boundary data
$h=u\chr{\Gs_\gb(A)}$ on $\Gs_\gb$.
\end{description}
\par Recall
that a  solution $u$ is moderate if $\abs{u}$ is dominated by a
harmonic function. When this is the case, $u$ possesses a boundary trace  (denoted by $\tr u$) given
by a bounded Borel measure. The boundary trace is attained in the sense of weak convergence,
as in the case of positive harmonic functions   (see \cite {MV1} and the references therein).
If $\tr u$ happens to be \abc
relative to Hausdorff $(N-1)$-dimensional measure on $\bdw$ we refer to its density $f$ as the
$L^1$ boundary trace of $u$ and write $\tr u=f$ (which should be seen as an abbreviation for
$\tr u= fd\,\BBH_{N-1}$).
\par A positive solution $u$ is \gsmod if there exists an increasing \seq of moderate solutions $\set{u_n}$
\sth $u_n\uar u$. This notion was introduced by Dynkin and Kuznetsov \cite{DK3} (see also \cite{Ku} and \cite{Dbook1}).

\par If $\mu$ is a bounded Borel measure on $\bdw$, the  problem
\begin{equation}\label{BVP}
 -\Gd u+u^q=0 \text{ in }\Gw, \q u=\mu \text{ on }\bdw
\end{equation}
possesses a (unique) solution \ifif $\mu$ vanishes on sets of
$\Capq$-capacity zero, (see \cite {MV3} and the references
therein). The solution is denoted by $u_\mu$.
\par
 The set of positive solutions of \rqq in $\Gw$ will be denoted by $\CU(\Gw)$. It is well
 known that this set is compact in the topology of $C(\Gw)$, i.e.,
 relative to local uniform convergence in $\Gw$.

 Our first result displays a dichotomy which is the basis for our definition of boundary trace.
\bth{dichotomy}
 Let $u$ be a positive solution of \rqq and let $\gx\in \bdw$. Then,\1
 {\bf either,} for every $\GTT_q$-open \ngh $Q$ of $\,\gx$, we have
 \begin{equation}\label{dich1}
  \lim_{\gb\to 0}\int_{\Gs_\gb(Q)}udS=\infty
\end{equation}
{\bf or} there exists  a $\GTT_q$-open \ngh $Q$ of $\,\gx$ \sth
 \begin{equation}\label{dich2}
  \lim_{\gb\to 0}\int_{\Gs_\gb(Q)}udS<\infty.
\end{equation}
\bcom
\begin{description}
  \item[either] for every $\GTT_q$-open \ngh $Q$ of $\,\gx$,
  $\set{\ugb{Q}}$ is \mdiv as $\gb\to 0$,
\item[or] there exists  a $\GTT_q$-open \ngh $Q$ of $\,\gx$ \sth
$\set{\ugb{Q}}$ is \mcon as $\gb\to0$.
\end{description}
\end{comment}
The first case  occurs
\ifif
\begin{equation}\label{int-sing}
\int_{D}u^q\gr(x)dx=\infty,\q D=(0,\gb_0)\ti Q
\end{equation}
for every $\GTT_q$-open \ngh $Q$ of $\gx$.
\es

A point $\gx\in\bdw$ is called a {\em singular} point of $u$ in
the first case, i.e. when  \req{dich1} holds, and
 a {\em regular} point of $u$ in the second case.
The set of singular points is denoted by $\CS(u)$ and its
complement in $\bdw$ by $\CR(u)$.
\bcom

\begin{description}
  \item[a]
 If $u$ is a solution of \rqq in $\Gw$ and $Q\sbs \bdw$ is a $\GTT_q$ open set we denote by $u_\gb^Q$, $0<\gb<\gb_0$,
 the solution of \rqq in $\Gw'_\gb$ with boundary data $h(x)=u(x)\chr{Q}(\gs_x)$ on $\Gs_\gb$.
 If $Q=\bdw$ the upper index will be omitted.
\item[b] The limit set of $\set{u_\gb^Q}$ as $\gb\to0$ is
 the set of functions $w$ in $\Gw$ \sth $u_{\gb_n}^Q\to w$ (locally uniformly in $\Gw$) for some \seq of positive numbers
 $\set{\gb_n}$ converging to zero. This limit set will be denoted by $\Lim_{\gb\to 0}u_\gb^Q$.
 If $u\in \CU(\Gw)$, the compactness of $\CU(\Gw)$ implies that, $\Lim_{\gb\to 0}u_\gb^Q$ is never empty and the set is a singleton
 \ifif $\lim_{\gb\to 0}u_\gb^Q$ exists.
\item[c] With $u$ and $Q$ as above we denote:
\begin{equation}\label{[u]^Q}
 [u]_0^Q:=\sup \set{w:\, w\in\Lim_{\gb\to 0}u^A_\gb, \q A \text{ $\GTT_q$-open,}\q\tl A\sbs Q}
\end{equation}
The supremum of such a set is a solution of \rqq.
\end{description}
\end{comment}
\par Our next result provides additional information on the behavior of solutions near the regular boundary set $\CR(u)$.

\bth{int-B} The set of regular points $\CR(u)$
is $\GTT_q$-open and there exists a non-negative Borel measure
$\mu$ on $\bdw$ possessing the following properties.\1
(i) For every $\gs\in \CR(u)$ there exist
a $\GTT_q$-open \ngh $Q$ of $\gs$ and a moderate solution $w$ \sth
\begin{equation}\label{gs-ngh}
\tl Q\sbs \CR(u), \q \mu(\tl Q)<\infty,
\end{equation}
and
 \begin{equation}\label{intlim}
 u^{Q}_\gb\to w \txt{locally uniformly in $\Gw$}, \q (\tr w)\chr{Q}= \mu\chr{Q}.
\end{equation}
(ii) $\mu$ is outer regular relative to $\GTT_q$.
\es
Based on these results we define the {\em precise boundary
trace} of $u$ by
\begin{equation}\label{couple-tr}
\tr^c u=(\mu,\CS(u)).
\end{equation}
 Thus a trace is represented by a couple
$(\mu,\CS)$, where $\CS\sbs \bdw$ is $\GTT_q$-closed and $\mu$ is an outer regular measure
relative to $\GTT_q$ which is
$\GTT_q$-locally finite on $\CR=\bdw\sms \CS$.
 However, not every couple of this type is a trace.
A necessary and sufficient condition for such a couple to be a trace is provided in
\rth{existence}.

\par The trace can also be represented by a Borel measure
$ \nu$ defined as follows:
\begin{equation}\label{meas-tr}
\nu(A)=
  \begin{cases}
    \mu(A) & \text{if }A\sbs \CR(u), \\
    \infty& \text{otherwise},
  \end{cases}
\end{equation}
for every Borel set $A\sbs \bdw$. We put
\begin{equation}\label{tr_u}
 \tr u:=\nu.
\end{equation}
This measure has the following properties:\1
(i) It is outer regular relative to $\GTT_q$. \\
(ii) It is {\em essentially absolutely continuous} relative to $\Capq$, i.e.,
for every $\GTT_q$-open set $Q$ and every Borel set $A$ \sth $\Capq(A)=0$,
 $$\nu(Q)=\nu(Q\sms A).$$

The second property will be denoted by $\nu\ppf\Capq$. It implies that, if $\nu(Q\sms A)<\infty$
then $\nu(Q\cap A)=0$. In particular, $\nu$ is \abc relative to $\Capq$ on $\GTT_q$-open sets
on which it is bounded.
\par A positive Borel measure possessing properties (i) and (ii) will be called a {\em $q$-perfect measure}.
The space of $q$-perfect measures will be denoted by $\BBM_q(\bdw)$.

\par We have the following necessary and sufficient condition for existence:
\bth{exist} Let $\nu$ be a positive Borel measure on $\bdw$,
possibly unbounded.
\par The boundary value problem
\begin{equation}\label{int-BVP}
 -\Gd u + u^q=0,\q u>0 \text{ in } \Gw,\q \tr (u)=\nu \text{ on }\bdw
\end{equation}
possesses a solution \ifif $\nu$ is $q$-perfect. When this
condition holds, a solution of \req{int-BVP} is given by
\begin{equation}\label{precisely}
 U=v\oplus U_{F}, \q v=\sup\set{u_{\nu\chr{Q}}:Q\in \CF_\nu},
\end{equation}
where
$$\CF_\nu:=\set{Q : \, Q\textrm{ \qop }, \;\nu(Q)<\infty},\q G:=\bigcup_{\CF_nu}Q,\q F=\bdw\sms G$$
and $U_F$ is the maximal solution vanishing on $\bdw\sms F$. \es

Finally we establish the following uniqueness result.
\bth{unique} Let $\nu$ be a  $q$-perfect measure on $\bdw$. Then the solution $U$ of problem
\req{int-BVP} defined by \req{precisely}
is \gsmod and it is the maximal solution with boundary trace $\nu$.
Furthermore, the solution of \req{int-BVP} is unique in the class of \gsmod solutions.
\es

\par For $q_c\leq q\leq 2$,  results similar to those stated in the last two theorems,
were obtained by Dynkin and Kuznetsov \cite{DK3} and Kuznetsov \cite{Ku},
based on their definition of fine trace.
  However, by their results, the prescribed trace is attained
  only up to equivalence,  i.e., up to
a set of capacity zero. By the present results, the solution
attains precisely the prescribed trace and this holds for all
values of $q$ in the supercritical range. The relation between the
Dynkin-Kuznetsov definition (which is used in a probabilistic
formulation) and the definition presented here, is not yet clear.
\par The plan of the paper is as follows:\1
{\em Section 2} presents results on the $\Capq$-fine topology which, for brevity, is called
the $q$-topology.\\
{\em Section 3} deals with the concept of maximal solutions which vanish on the boundary
outside a
\qcl set. Included here is a sharp estimate for these solutions, based on the capacitary estimates developed
by the authors in \cite{MV4}. In particular we prove that the maximal solutions are \gsmod. This was established
in \cite{MV4} for solutions vanishing on the boundary outside a {\em compact set}.\\
{\em Section 4}  is devoted to the problem of localization of solutions in terms of boundary
behavior. Localization methods are of crucial importance in the study of trace and the associated boundary
value problems. The development of these methods is particularly subtle in the supercritical case.\\
{\em Section 5} presents the concept of precise trace and studies it, firstly on the regular boundary set, secondly
in the case of \gsmod solutions and finally in the general case. This section contains the proofs of the
theorems stated above:
\rth{dichotomy} is a consequence of \rth{loc-tr2}.
\rth{int-B} is a consequence of \rth{reg-tr}. \rth{exist} is a consequence of \rth{existence} (see the remark
following the proof of the latter theorem). Finally \rth{unique}
is contained in \rth{existence}.

\section{The $q$-fine topology}
A basic ingredient in our study is the fine topology
associated with a Bessel capacity on $(N-1)$-dimensional smooth manifolds. 
 The theory of fine topology associated  with the Bessel capacity
$C_{\ga,p}$ in $\BBR^N$ essentially requires $0<\ga p\leq N$ (see
\cite[Chapter 6]{AH}). In this paper we are interested in the fine
topology associated  with the capacity $\Capq$ in $\BBR^{N-1}$ or
on the boundary manifold $\bdw$ of a smooth bounded domain $\Gw
\in \BBR^N$. We assume that $q$ is in the supercritical range for
\rqq, i.e., $q\geq q_c=(N+1)/(N-1).$ Thus $2q'/q=2/(q-1)\leq N-1$.
We shall refer to the $(2/q,q')$-fine topology briefly as the
$q$-topology.

 An important concept related to this topology is the
$(2/q,q')$-quasi topology. We shall refer to it as the $q$-quasi
topology. For definition and details see \cite[Section~6.1-4]{AH}.

We say that a subset of $\bdw$  is $q$-open (resp. $q$-closed) if
it is open (resp. closed) in the $q$-topology on $\bdw$. The terms
$q$-quasi open and $q$-quasi closed are understood in an analogous
manner.\2 \note{Notation 2.1} Let $A,B$ be subsets of $\BBR^{N-1}$
or of $\bdw$.
\begin{description}
\item[{\rm a}] $A$ is $q$-essentially contained in $B$, denoted  $A \overset{q}{\sbs} B$, if
$$\Capq(A\sms B)=0.$$
\item[{\rm b}] The sets $A, B$ are $q$-equivalent, denoted $A \overset{q}{\sim} B$, if
$$\Capq(A\Gd B)=0.$$
\item[{\rm c}]  The $q$-fine closure of a set $A$ is denoted by $\tl A$. The $q$-fine interior of $A$
is denoted by $A^\diamond$.
\item[{\rm d}] Given $\ge>0$, $A^\ge$ denotes the intersection of $\BBR^{N-1}$ (or $\bdw$) with the $\ge$-\ngh  of
$A$ in $\BBR^N$.
\item[{\rm e}] The set of
$(2/q,q')$-thick (or briefly $q$-thick) points of $A$ is denoted
by $b_q(A)$. The set of $(2/q,q')$-thin (or briefly $q$-thin)
points of $A$ is denoted by $e_q(A)$, (for  definition see
\cite[Def. 6.3.7]{AH}).
\end{description}
\Remark If $A\sbs \bdw$  and $B:=\bdw\sms A$ then
\begin{equation}\label{eqB}
 A \textrm{ is \qop }\Llra A\sbs e_q(B),\q B \textrm{ is \qcl }\Llra b_q(B)\sbs B.
\end{equation}
\Consy
\begin{equation}\label{tlA}
 \tl A=A\cup b_q(A),\q  A^\diamond=A\cap e_q(\bdw\sms A),
\end{equation}
(see \cite[Section 6.4]{AH}.)\2
\indent  The  capacity $\Capq$ possesses the Kellogg property, namely,
\begin{equation}\label{A-b(A)}
 \Capq(A\sms b_q(A))=0,
\end{equation}
(see \cite[Cor. 6.3.17]{AH}). Therefore
\begin{equation}\label{b_qA}
A\sbsq b_q(A)\smq \tl A
\end{equation}
but, in general, $b_q(A)$ does not contain $A$. The Kellog
property and \req{eqB} implies: \bprop{qop/qcl}
(i) If $Q$ is a \qop set then $\check Q:=e_q(\bdw\sms Q)$ is the largest \qop set
that is $q$-equivalent to $Q$.\\
(ii) If $F$ is a \qcl set then $\hat F=b_q(F)$ is the smallest
\qcl set that is $q$-equivalent to $F$.
\es

We collect below
several facts concerning the $q$-fine topology that are used
throughout the paper. \bprop{q-top} Let $q_c=(N+1)/(N-1)\leq q$.
\begin{description}
\item[{\rm i}] Every $q$-closed set is $q$-quasi closed \cite[Prop.~6.4.13]{AH}.
\item[{\rm ii}] If $E$ is $q$-quasi closed then  $\tl E\smq E$ \cite[Prop.~6.4.12]{AH}.
\item[{\rm iii}] A set $E$ is $q$-quasi closed \ifif there exists a \seq  $\{E_m\}$ of
compact subsets of $E$ \sth $\Capq(E\sms E_m)\to 0$ \cite[Prop.
6.4.9]{AH}.
\item[{\rm iv}] There exists a constant $c$ \sth, for every set $E$,
$$\Capq(\tl E)\leq c \Capq (E),$$
see \cite[Prop. 6.4.11]{AH}.
\item[{\rm v}] If $E$ is $q$-quasi closed and $F\smq E$ then $F$ is $q$-quasi closed.
\item[{\rm vi}] If $\set{E_i}$ is an increasing \seq of arbitrary sets then
$$\Capq(\cup E_i)=\lim_{i\tin}\Capq(E_i).$$
\item[{\rm vii}] If $\set{K_i}$ is a decreasing \seq of compact sets then
$$\Capq(\cap K_i)=\lim_{i\tin}\Capq(K_i).$$
\item[{\rm viii}] Every Suslin set and, in particular, every Borel set $E$ satisfies
$$\begin{aligned}
\Capq(E)=\,&\sup\set{\Capq(K):\, K\sbs E,\; K\text{ compact}}\\
=\,&\inf\set{\Capq(G):\, E\sbs G,\; G\text{ open}}.
\end{aligned}$$
\end{description}
\es

For the last three statements see \cite[Sec. 2.3]{AH}. Statement (v)
 is an easy consequence of \cite[Prop. 6.4.9]{AH}. However note that this assertion
is no longer valid if '$q$-quasi closed' is replaced by '\qcl'. Only the following weaker
statement holds:\1
\indent{\em If $E$ is \qcl and $A$ is a set \sth $\Capq(A)=0$ then $E\cup A$ is \qcl.}


\bdef{e-ngh}{\rm Let $E$ be a quasi closed set. An increasing \seq
$\{E_m\}$ of compact subsets of $E$ \sth $\Capq(E\sms E_m)\to 0$
is called a  $q$-{\em stratification of} $E$.\\
(i)  We say that  $\{E_m\}$ is
 a {\em proper $q$-stratification of} $E$ if
 \begin{equation}\label{strat}
  C_{2/q,q'}(E_{m+1}\sms E_m)\leq 2^{-m-1} C_{2/q,q'}(E).
\end{equation}
(ii) Let $\{\ge_m\}$ be a strictly decreasing \seq of positive
numbers converging to zero \sth
\begin{equation}\label{ngh-1}
 \Capq(G_{m+1}\sms G_m)\leq 2^{-m} C_{2/q,q'}(E),   \q
 G_m:=\cup_{k=1}^m E_k^{\ge_k}.
\end{equation}
The \seq $\set{\ge_m}$ is called a $q$-{\em proper \seq.}\\
(iii) If $V$ is a \qop set \sth $\Capq(E\sms V)=0$ we say that $V$ is a $q$-{\em quasi \ngh} of $E$.
}\es
\Remark Observe that $G:=\cup_{k=1}^\infty E_k^{\ge_k}$ is a \qop \ngh of
$E'=\cup E_m$ but, in general,  only a $q$-quasi \ngh of $E$.
\blemma{e-ngh} Let $E$ be  a $q$-closed set \sth $\Capq(E)>0$. Then: \1
(i) Let $D$ be an open set \sth $\Capq(E\sms D)=0$. Then $E\cap D$ is $q$-quasi closed
and \consy there exists a proper $q$-stratification of $E\cap D$, say $\{E_m\}$.
Furthermore, there exists a
$q$-proper \seq $\{\ge_m\}$  \sth
$$G=\cup_{m=1}^\infty (E_m)^{\ge_m}\sbs D$$
and
\begin{equation}\label{ngh-2a}
 \cup E_m=E'\sbs O\sbs \tl O \sbsq D \txt{where} O:=\cup_{m=1}^\infty (E_m)^{\ge_m/2}.
\end{equation}
\Consy
\begin{equation}\label{ngh-2}
 E\sbsq O\sbs \tl O \sbsq D.
\end{equation}
(ii) If $D$ is a \qop set \sth $E\sbsq D$ then there exists a \qop
set $O$ \sth \req{ngh-2} holds.
\es
\bcom (iii) If $K$ is a
non-empty compact set \sth $K\sbs D$ then there exists a \qop set
$O$ \sth
\begin{equation}\label{ngh-2'}
 K\sbs O\sbs\tl O \sbsq D.
\end{equation}
 \Remark Note that in (iii) we do not exclude the case $\Capq(K)=0$ but we require $K\sbs D$, not just $K\sbsq D$,
and the conclusion is stronger then in part (ii): $K\sbs O$.
\end{comment}
 \bproof  If $A_1,A_2$ are two sets \sth $A_1\smq A_2$ and $A_1$ is $q$-quasi closed then $A_2$ is $q$-quasi closed,
 (see the discussion of the quasi topology in \cite[sec. 6.4]{AH}).
 Since $E\cap D\smq E$ and $E$ is  \qcl it follows  that
 $E\cap D$ is $q$-quasi closed. Let $\{E_m\}$ be a proper $q$-stratification of $E\cap D$ and put
 $E'=\cup_{m=1}^\infty E_m$. If $E'$ is a closed set the remaining part of assertion (i) is trivial.
 Therefore we assume that $E'$ is not closed and that $$\Capq(E_{m+1}\sms E_m)>0.$$
To prove the first statement we  construct the sequence
$\{\ge_m\}_{m=1}^\infty$  inductively
 so that (with $E_0=\ems$ and $\ge_0=1$) the following conditions are satisfied:
\begin{align}\label{ngh-3}
& F_m:=E_m\sms (E_{m-1})^{\rec{2}\ge_{m-1}}, &&  \Capq(F_m)>0,  \\
&  \Capq( \overline{F_m^{\ge_m}})\leq 2\Capq(E_m\sms E_{m-1}),  && \overline{F_m^{\ge_m}}\sbs D, \label{ngh-4}\\
  &  \ge_m<\ge_{m-1}/2, && m=1,2,\dots.\label{ngh-5}
\end{align}
Choose $0<\ge_1<1/2$, sufficiently small so that
$$E_1^{\ge_1}\sbs D, \q \Capq(E_2\sms E_1^{\ge_1/2})>0.$$
 This is possible because our assumption
implies that there exists a compact subset of $E_2\sms E_1$ of
positive capacity.  By induction we obtain
\begin{equation}\label{ngh-6}
 E_{m}^{\ge_{m}}\sbs E_{m-1}^{\ge_{m-1}}\cup  F_m^{\ge_m} \sbs D
\end{equation}
and consequently
\begin{equation}\label{ngh-7}
 E_m^{\ge_m}\sbs\cup_{k=1}^m F_k^{\ge_k}, \q m=1,2,\dots\,.
\end{equation}
 Since $F_m\sbs E_m$, \req{ngh-7} implies that
\begin{equation}\label{ngh-8}
 G:=\cup_{k=1}^\infty E_k^{\ge_k}=\cup_{k=1}^\infty F_k^{\ge_k},\q
 G_m:=\cup_{k=1}^m E_k^{\ge_k}=\cup_{k=1}^m F_k^{\ge_k}.
\end{equation}
The \seq $\{\ge_m\}$ constructed above satisfies \ref{ngh-1}.
Indeed, by \req{strat},  \req{ngh-4} and \req{ngh-8},
\begin{align}\label{ngh-9}
 \Capq(G\sms G_m)&\leq \sum_{k=m+1}^\infty\Capq(F_k^{\ge_k})\\
 &\leq 2\sum_{k=m+1}^\infty\Capq(E_k\sms E_{k-1})\leq 2^{-m+1}\Capq(E).\notag
\end{align}
\par Next we show  that the set
$$O':=\cup_{k=1}^\infty \overline{E_k^{\ge_k/2}}$$
is $q$-quasiclosed. By \req{ngh-7},
\begin{equation}\label{ngh-10}
O'= \cup_{k=1}^\infty \overline{F_k^{\ge_k/2}}, \q
O'_m:=\cup_{k=1}^m \overline{E_k^{\ge_k/2}}=\cup_{k=1}^m
\overline{F_k^{\ge_k/2}}.
\end{equation}
Hence, by \req{strat} and \req{ngh-4},
\begin{equation}\label{ngh-11}
\Capq(O'\sms O'_m)\leq
\sum_{k=m+1}^\infty\Capq(\overline{F_k^{\ge_k/2}})\leq
2^{-m+1}\Capq(E).
\end{equation}
Since $O'_m$ is closed this implies that $O'$ is quasiclosed.
 Further any quasiclosed set is equivalent to its fine closure.
 Since $O\subset O'$ it follows that $\tl O\sbs \tl O'\smq O'\sbs G$.
 \par We turn to the proof of  (ii)  for which we need the following:\2
 \note{Assertion 1} Let $D$ be a \qop set. Then there exists a \seq  of relatively open sets $\set{A_n}$   \sth
\begin{equation}\label{A_n}
D_n:=D\cup A_n \text{ is open},\q \Capq(\tl A_n)\leq 2^{-n}, \q
\tl A_{n+1}\sbsq A_n.
\end{equation}
\par The \seq is constructed inductively. Let $D'_1$ be an open set \sth $D\sbs D'_1$ and $A'_1=D'_1\sms D$ satisfies
$\Capq (\tl A'_1)<1/4$. Let $A_1$ be an open set \sth $\tl
A'_1\sbs A_1$ and $\Capq(\tl A_1)\leq 1/2$. Assume that we
constructed $\set{A'_k}_1^{n-1}$ and $\set{A_k}_1^{n-1}$ so that
the sets $A_k$ are open, \req{A_n} holds and
\begin{equation}\label{A'_n}\BAL
D'_k=D\cup A'_k \text{ is open}, &\q \set{D'_k}_1^{n-1} \text{ is decreasing}, \\
 \Capq (\tl A'_{k})<2^{-(k+1)}, &\q
\tl A'_k\sbs A_k. \EAL\end{equation} Let $D'_{n}$ be an open set
\sth
$$D\sbs D'_{n}\sbs D'_{n-1}, \q \Capq (\tl A'_{n})<2^{-(n+1)}\txt{where} A'_{n}=D'_{n}\sms D.$$
Then $A'_n\sbs A'_{n-1}$ and \consy $\tl A'_n\sbs A_{n-1}$.  Since
$A_{n-1}$ is open, statement (i) implies that
 there exists an open set $A_n$ \sth
$$\tl A'_n\sbsq A_n\sbs \tl A_n\sbsq A_{n-1}, \q \Capq(\tl A_n)\leq 2^{-n}.$$
This completes the proof of the assertion.\2
 \indent  Let $A_n$ and $D_n$ be as in \req{A_n}. By (i) there exists a \qop set $Q$ \sth
$$E\sbsq Q\sbs \tl Q\sbsq D_1.$$
Put
$$Q_n:= Q\sms (\cup_1^{n-1}(\tl A_k\sms  A_{k+1}),\q Q_\infty=Q\sms (\cup_1^{\infty}(\tl A_k\sms  A_{k+1}).$$
Then $Q_n$ is a \qop set and we claim that
\begin{description}
\item[{\rm\bf a}]\hspace{5mm}$Q_\infty \text{ is quasi open,}$
 \item[{\rm\bf b}]\hspace{5mm}$E\sbsq Q_n\sbsq D_n,$
 \item[{\rm\bf c}]\hspace{5mm}$ \tl Q_n\sbsq D_n\cup (\cup_1^n\prt_qA_i)$
 \end{description}
\indent Since $Q_n$ is \qop, {\bf a} follows from \req{A_n} which
implies:
$$Q_\infty\cup (\cup_n^\infty(\tl A_k\sms  A_{k+1})=Q_n,\q \Capq((\cup_n^\infty(\tl A_k\sms  A_{k+1})\leq2^{n-1}.$$
We verify {\bf b, c} by induction. Put $Q_1=Q$ so that {\bf b, c} hold for $n=1$.
If {\bf b} holds for $n=1,\cdots,j$ then,
$$ Q_{j+1}= Q_j\sms  (\tl A_j\sms  A_{j+1}),\q E\sbsq Q_j\sbsq D_j,$$
which implies {\bf b} for $n=j+1$. If {\bf c} holds for $n=1,\cdots,j$
then,
$$\tl Q_{j+1}\sbs \tl Q_j\sms  ( A_j\sms  \tl A_{j+1})\sbsq D\cup (\cup_1^j\prt_qA_i)\cup \tl A_{j+1})=
D_{j+1}\cup (\cup_1^{j+1}\prt_qA_i),$$ so that {\bf c} holds for
$n=j+1$.
\par Taking the limit in {\bf b} as $n\tin$ we obtain $$E\sbsq Q_\infty\sbsq D.$$
Taking the limit in {\bf c} as $n\tin$ we obtain
$$\tl Q_\infty\sbsq D\cup (\cup_1^\infty\prt_qA_i).$$
However, by the same token,
$$\tl Q_\infty\sbsq D\cup (\cup_k^\infty\prt_qA_i) \forevery k\in\BBN.$$
Therefore, by \req{A_n}, $\tl Q_\infty\sbsq D$. Thus (ii) holds
with $O=Q_\infty$.
\eproof
\blemma{q-compact} Let $E$ be a \qcl set and let $\CD$ be a cover
of $E$ consisting of \qop sets. Then, for every $\ge>0$ there
exists an open set $O_\ge$ \sth $\Capq(O_\ge)<\ge$ and $E\sms
O_\ge$ is covered by  a finite subfamily of $\CD$. \es \bproof By
\cite[Sec. 6.5.11]{AH} the $(\ga,p)$-fine topology possesses the
quasi Lindel\"{o}f property. Thus there exists a denumerable
subfamily of $\CD$, say $\set{D_n}$, \sth
$$O=\cup\set{D: D\in \CD}\smq \cup D_n.$$
Let $O_n$ be an open set containing $D_n$ \sth $\Capq(O_n\sms
D_n)<\ge/(2^n3)$. Let $K$ be a compact subset of $E\cap(\cup_1^\infty D_n)$ \sth
$\Capq(E\sms K)<\ge/3$. Then $\set{O_n}$ is an open cover of $K$
so that there exists a finite subcover of $K$, say
$\set{O_1,\cdots, O_k}$. It follows that
$$\Capq(E\sms \cup_{n=1}^k D_n)\leq \Capq(E\sms K)+ \sum_n \Capq(O_n\sms D_n)<2\ge/3.$$
Let $O_\ge$ be an open subset of $\bdw$ \sth $E\sms \cup_{n=1}^k D_n\sbs O_\ge$ and $\Capq(O_\ge)<\ge$.
This set has the properties stated in the lemma.
\eproof
\blemma{set-approx} (a) Let $E$ be a q-quasi closed set
and $\{E_m\}$
 a proper $q$-stratification for $E$. Then  there exists a {\em decreasing} \seq of
open sets $\{Q_j\}$ \sth $\cup E_m:=E' \sbs Q_j$ for every $j\in
\BBN$ and
\begin{equation}\label{q-closed-app}\BAL
&(i)\q \cap Q_j= E', \q \tl Q_{j+1}\sbsq Q_j,\\
&(ii)\q  \lim \Capq(Q_j)= \Capq(E).
\EAL\end{equation} (b) If $A$ is a $q$-open set, there exists a
decreasing \seq of open sets $\set{A_m}$ \sth
\begin{equation}\label{open-strat}
 A\sbs \cap A_m=:A', \qq \Capq(A_m\sms A')\to 0, \qq A\smq A'.
\end{equation}
Furthermore there exists an {\em increasing} \seq of closed sets
$\set{F_j}$ \sth $F_j\sbs A'$ and
\begin{equation}\label{q-open-app}\BAL
&(i)\q \cup F_j= A', \q F_j\sbsq F_{j+1}^\diamond=:D_{j+1}, \\
&(ii)\q \Capq(F_j)\to \Capq(A)\phantom{=\Capq(\tl Q_j)}.
\EAL\end{equation} \es
\bproof (a)  Let  $\{\ge_m\}$ be a \seq of
positive numbers decreasing to zero satisfying \req{ngh-1}. Put
$$Q_j:=\cup_{m=1}^\infty (E_m)^{\ge_m/2^j}.$$
Then $E'=\cap Q_j$ and
\begin{equation}\label{Q'j}
 \tl Q_j\sbsq Q'_j:=\cup_{m=1}^\infty \ovl{E_m^{\ge_m/2^j}}\sbs Q_{j-1}.
\end{equation}
 Indeed $Q'_j$ is quasi closed so that  $\tl Q'_j\smq Q'_j$.
This proves \req{q-closed-app}(i).
\par If $D$ is a \ngh of $E'$ then, for every $k$ there exists $j_k$ \sth
$$\cup_{m=1}^k (E_m)^{\ge_m/2^j}\sbs D \forevery j\geq j_k.$$
Therefore,
$$\Capq(Q_j\sms D)\leq\Capq(Q_j\sms \cup_{m=1}^k (E_m)^{\ge_m/2^j})\leq 2^{-k+1}\Capq(E) \forevery j\geq j_k.$$
Hence
\begin{equation}\label{QjD}
  \Capq(Q_j\sms D)\to 0 \txt{as} j\tin.
\end{equation}
Let $\{D_i\}$ be a decreasing \seq of open neighborhoods of $E'$
\sth
$$\Capq(D_i)\to \Capq(E').$$
 By \req{QjD}, for every $i$ there exists $j(i)>i$ such that
\begin{equation}\label{QiDi}
  \Capq(Q_{j(i)}\sms D_i)\to 0 \txt{as} i\tin.
\end{equation}
 It follows that
 $$\Capq(E')\leq \lim \Capq(Q_{j(i)})\leq\lim\Capq(D_i)= \Capq(E')=\Capq(E).$$
  This proves \req{q-closed-app} (ii).\1
(b) Put $E=\bdw\sms A$ and let $\{E_m\}$ and $\{\ge_m\}$ be as in
(a). Then \req{open-strat} holds with $A_m:=\bdw\sms E_m$. In
addition, \req{q-open-app}(i) with $F_j:= \bdw\sms Q_j$ is a
consequence of \req{q-closed-app}(i).

To verify \req{q-open-app} (ii) we observe that, if $K$ is a
compact subset of $A'$ then, by \req{QjD},
$$\Capq(K\sms F_j)\to 0.$$
Let $\{K_i\}$ be an increasing \seq of compact subsets of $A'$
\sth
$$\Capq(K_i)\uparrow \Capq(A')=\Capq(A).$$
As in part (a), for every $i$ there exists $j(i)>i$ such that
\begin{equation}\label{KiFi}
  \Capq(K_i\sms F_{j(i)})\to 0 \txt{as} i\tin.
\end{equation}
 It follows that
 $$\Capq(A')\geq \lim \Capq(F_{j(i)})\geq \lim\Capq(K_i)= \Capq(A')=\Capq(A).$$
This proves \req{q-open-app} (ii). \eproof
\blemma{closedngh} Let $Q$ be a \qop set. Then, for every $\gx\in
Q$,
 there exists a \qop set $Q_\gx$ \sth
$$\gx\in Q_\gx\sbs \tl Q_\gx\sbs Q.$$
\es \bproof By definition, every point in $Q$ is a $q$-thin point
of $E_0=\bdw\sms Q$. Assume that $\diam Q<1$ and put:
$$r_n=2^{-n}, \q K_n=\set{\gs:r_{n+1}\leq\abs{\gs-\gx}\leq r_n},\q E_n:=E_0\cap K_n\cap \bar B_1(\gx).$$
Thus $E_n$ is a \qcl set; we denote $E=\cup_{n=0}^\infty E_n$. Since $\gx$ is a $q$-thin point of $E$,
$$\sum_0^\infty (r_n^{-N+1+2/(q-1)}\Capq(B_n\cap E))^{q-1}<\infty, \q B_n=B_{r_n}(\gx),$$
which is equivalent to
$$\sum_0^\infty \big(r_n^{-N+1+2/(q-1)}\Capq(E_n)\big)^{q-1}<\infty.$$

\par Let $\set{E_{m,n}}_{m=1}^\infty$ be a $q$-proper stratification of $E_n$.  
Let $\bar\ge^n:=\set{\ge_{m,n}}_{m=1}^\infty$ be a $q$-proper \seq
(relative to the above stratification) \sth $\ge_{1,n}\in (0,r_{n+2})$ and
$$\Capq( V_n)<2\Capq(\cup E_n) \txt{where} V_n:=\cup_{m=1}^\infty E_{m,n}^{\ge_{m,n}}\cap B_1(\gx).$$
Then $V_n\sbs K_{n-2}\sms K_{n+2}$,
 $\gx$ is a $q$-thin point of the set $G=\cup_0^\infty V_n$ and $\gx\nin G$.
\Consy $\gx\nin\tl G$.
\par Put
$$Z_n:=\cup_{m=1}^\infty \overline{ E_{m,n}^{\ge_{m,n}/2}}\cap\bar B_{1/2}(\gx),\q
F_0:=\cup _0^\infty  Z_n.$$
Since $Z_n\sbs V_n$ it follows that $\gx$ is a $q$-thin point of $F_0$ and $\gx\nin \tl F_0$.
\Consy $Q_0:=(Q\cap B_{1/2}(\gx))\sms \tl F_0$ is a \qop subset of $Q$ \sth
$$\gx\in Q_0, \q \tl Q_0\sbs (\tl Q\cap\bar B_{1/2}(\gx)) \sms F_0\sbs
(\tl Q\cap\bar B_{1/2}(\gx))\sms E\sbs Q.$$
\eproof
\section{Maximal solutions}

We consider positive solutions of the equation \rqq with $q\geq
q_c$, in a  bounded domain $\Gw\sbs \BBR^N$ of class $C^2$. A
function $u\in L^q_{loc}(\Gw)$ is a subsolution  (resp.
supersolution) of the equation if $-\Gd u+|u|^{q-1}u\leq 0$ (resp.
$\geq 0$) in the distribution sense.

If $u\in L^q_{loc}(\Gw)$ is a subsolution of the equation then (by
Kato's inequality \cite{Ka}) $\Gd|u|\geq |u|^q$. Thus $|u|$ is
subharmonic and  \consy $u\in L^\infty_{loc}(\Gw)$. If $u\in
L^q_{loc}(\Gw)$
 is a solution then $u\in C^2(\Gw)$.

 An increasing  \seq of bounded domains of class $C^2$, $\set{\Gw_n}$,  \sth $\Gw_n\uparrow \Gw$ and
$\bar \Gw_n\sbs \Gw_{n+1}$ is called an {\em exhaustive} \seq
relative to $\Gw$.

\bprop{sub-sup} Let $u$ be a non-negative function in
$L^\infty_{loc}(\Gw)$.\1 (i) If $u$ is a subsolution of
\req{q-eq}, there exists a minimal solution
 $v$ dominating $u$, i.e.,  $u\leq v\leq U$ for any solution $U\geq u$.\1
(ii) If $u$ is a supersolution of \req{q-eq}, there exists a
maximal solution $w$ dominated by $u$, i.e., $V<w<u$  for any
solution $V\leq u$.

All the inequalities above are a.e.\,.
\es
\bproof
Let $u_\ge=J_\ge u$ where $J_\ge$ is a smoothing operator and $u$
is extended by zero outside $\Gw$. Put $\tl u=\lim_{\ge\to 0}
u_\ge$ (the limit exists a.e. in $\Gw$ and $\tl u=u$ a.e.).
Let $\gb_0$, $\Gw_\gb$, $\Gs_\gb$ etc. be as in  Notation 1.1.
Since $u_\ge\to \tl u$ in $L^1(\Gw)$ it follows that
$$u_\ge\rest{\Gs_\gb}\to \tl u\rest{\Gs_\gb} \txt{in} L^1(\Gs_\gb)$$
for a.e.$\gb\in (0,\gb_0)$. Choose a \seq $\set{\gb_n}$ decreasing to zero \sth the above
convergence holds for each surface $\Gs_n:=\Gs_{\gb_n}$. Put $D_n:=\Gw'_{\gb_n}$.
Assuming that
$u$ is a \subsol of \rqq in $\Gw$, $u_\ge$ is a subsolution of the
boundary value problem for \rqq in $D_n$ with boundary data
$u_\ge\rest{\Gs_n}$. \Consy $\tl u$ is a subsolution of the
boundary value problem for \rqq in $D_n$ with boundary data
$\tl u\rest{\Gs_n}\in L^1(\Gs_n)$. (Here we use the assumption
$u\in L^\infty_{loc}(\Gw)$ in order to ensure that $u_\ge^q\to \tl u^q$ in $L^1_{loc}(\Gw)$.)
Let $v_n$ denote the solution of this boundary value problem in the $L^1$ sense:
$$-\Gd v_n + v_n^q=0 \txt{in} D_n, \qq v_n=\tl u \txt{on}\Gs_n.$$
Then $v_n\in C^2(D_n)\cap L^\infty(D_n)$, $v_n\leq \norm{u}_{L^\infty(D_n)}$ and the boundary
data is assumed in the $L^1$ sense.
Clearly $\tl u\leq v_n$ in $D_n$, n=1,2,\dots\,. In particular, $v_{n}\leq v_{n+1}$ on $\Gs_n$.
 This implies $v_n\leq v_{n+1}$ in $D_n$. In addition, by the Keller-Osserman inequality
  the \seq $\set{v_n}$ is eventually bounded in every compact subset of $\Gw$.
Therefore $v=\lim v_n$ is the \sol with the properties  stated in (i).

Next assume that $u$ is a \supsol and let $\set {D_n}$ be as
above. Since $u\in L^q(D_n)$ there exists a positive solution
$w_n$ of the boundary value problem
$$-\Gd w=u^q \text{ in }D_n,\q w=0 \text{ on }\Gs_n.$$
Hence $u+w_n$ is superharmonic and  its boundary trace is
precisely $\tl u\rest{\Gs_n}$. \Consy $u+w_n\geq z_n$ where $z_n$
is the harmonic function in  $D_n$ with boundary data $\tl
u\rest{\Gs_n}$. Thus $u_n:= z_n-w_n$ is the smallest solution of
\rqq in $D_n$ dominating $u$. This implies that $\set{u_n}$
decreases and the limiting solution $U$ is the smallest solution
of \rqq dominating $u$ in $\Gw$. \eproof
\bprop{ex-sub} Let $u,v$ be non-negative, locally bounded
functions in $\Gw$.\1 (i) If $u,v$ are subsolutions (resp.
supersolutions) then $\max(u,v)$ is a \subsol (resp. $\min(u,v)$
is a \supsol).\1 (ii) If $u,v$ are  supersolutions then $u+v$ is a
\supsol.\1 (iii) If $u$ is a \subsol and $v$ a \supsol then
$(u-v)_+$ is a \subsol. \es \bproof The first two statements are
well known; they can be verified by an application of Kato's
inequality. The third statement is verified in a similar way:
$$\Gd(u-v)_+=\sign_+(u-v)\Gd(u-v)\geq (u^q-v^q)_+\geq (u-v)_+^q.$$
\eproof
\note{Notation 3.1} Let $u,v$ be non-negative, locally bounded
functions in $\Gw$.\1 (a) If $u$ is a \subsol, $[u]_\dag$ denotes
the smallest solution dominating $u$.\1 (b) If $u$ is a \supsol,
$[u]^\dag$ denotes the largest solution dominated by $u$.\1 (c) If
$u,v$ are subsolutions then $u\vee v:=[\max(u,v)]_\dag$.\1 (d) If
$u,v$ are supersolutions then $u\wedge v:=[\inf(u,v)]^\dag$ and
$u\oplus v:=[u+v]^\dag$.\1 (e) If $u$ is a \subsol and $v$ a
\supsol then $u\ominus v:=[(u-v)_+]_\dag$.\2
 The following result was proved in \cite{Ku} (see also
\cite[Sec.~8.5]{Dbook1}).
\bprop{sup u_n} (i) Let $\set{u_k}$ be a
\seq of positive, continuous subsolutions of \req{q-eq}. Then
$U:=\sup u_k$ is a subsolution. The statement remains valid if
\subsol is replaced by \supsol and $\sup$ by $\inf$.\1 (ii) Let
$\CT$ be a family of positive solutions of \req{q-eq}. Suppose
that, for every pair $u_1,u_2\in \CT$, there exists $v\in \CT$
\sth
$$\max(u_1,u_2)\leq v,\txt{resp.}  \min(u_1,u_2)\geq v.$$
Then there exists a monotone \seq $\set{u_n}$ in $\CT$ \sth
$$u_n\uparrow\sup\CT, \txt{resp.} u_n\downarrow\sup\CT.$$
Thus $\sup\CT$ (resp. $\inf \CT$) is a solution. \es \bcom \Remark
The result is due to \cite{Ku} (see also \cite[Ch.8,
Sec.5]{Dbook1}. For the convenience of the reader we provide a
brief proof. \bproof Let $D$ be a smooth domain \sth $D\Subset
\Gw$ and let $v_n$ be the solution of \req{q-eq} in $D$ \sth
$v_n=\max(u_1,\dots,u_n)$ on $\prt D$. Then $v=\lim v_n$ is a
solution in $D$ \sth $v=\sup u_n$ on $\prt D$ and it is
 the minimal solution which dominates $\sup u_n$ in $D$. This implies (i).

We turn to the second assertion. For every point $\gx\in \Gw$
there exists a \seq $\set{u_n}$ in $\CT$
 \sth $u_n(\gx)\uparrow U(\gx)$. The assumptions of (ii) imply that $\set{u_n}$ can be chosen to
 be an increasing sequence. \Consy $w_\gx=\sup u_n$ is a solution.
  Let $\set{\gx_k}$ be a dense countable subset of $\Gw$
 and let $w=\sup w_{\gx_k}$. Then $w$ is a \subsol dominated by $U$ \sth $w(\gx_k)=U(\gx_k)$ for every $k$.

Let $\set{u_n}$ be a sequence of positive subsolutions of
\req{q-eq} and put $v_n=\vee_{k=1}^n u_k$. Then $\set{v_n}$ is a
monotone increasing \seq of solutions so that $v=\lim v_n$ is a
solution of \req{q-eq}. \eproof
\end{comment}


\bdef{vanishing}{\rm A solution $u$ of \req{q-eq} vanishes on a
relatively open set $Q\sbs \bdw$ if $u\in C(\Gw\cup Q)$ and $u=0$
on $Q$. A  positive solution $u$ vanishes on a $q$-open set $A\sbs
\bdw$ if
$$u=\sup\set{v\in \CU(\Gw):\, v\leq u,\; v=0 \;\textrm{on some relatively open \ngh of $A$}}.$$
 When this is the case we write $u\app{A}0$.
}\es
\blemma{vanishing} Let $A$ be a \qop subset of $\bdw$ and $u_1,u_2\in \CU(\Gw)$.\\
(a) If both solutions vanish on $A$ then $u_1\vee u_2\app{A}0$. If
$u_2\app{A}0$  and $u_1\leq u_2$ then $u_1\app{A}0$.\1 (b) If
$u\in \CU(\Gw)$ and $u\app{A}0$ then there exists an increasing \seq of
solutions $\{u_n\}\sbs \CU(\Gw)$, each of which vanishes on a
relatively open \ngh of $A$ (which may depend on $n$) \sth
$u_n\uparrow u$.\1 (c) If  $A,A'$ are \qop sets,  $A\smq A'$ and
$u\app{A}0$ then $u\app{A'}0$. \es \bproof The first assertion
follows easily from the  definition. Thus the set of solutions
$\{v\}$ described in the definition is closed \wrto the binary
operator $\vee$. Therefore, by \rprop{sup u_n}, the supremum of
this set is the limit of an increasing \seq of elements of this
set.
\par The last statement is obvious.
\eproof
\bdef{maxsol}{\rm (a) Let $u\in \CU(\Gw)$ and let $A$ denote the union of
all \qop sets on which $u$ vanishes. Then $\bdw\sms A$ is called
the {\em fine boundary support} of $u$, to
be denoted by $\supp^q_{\bdw}u$.\\[1mm]
(b) For any Borel set $E$  we denote
$$U_E=\sup\set{u\in \CU(\Gw): u\app{\tl E^c}0,\; \tl E^c=\bdw\sms \tl E}.$$
Thus $U_E=U_{\tl E}$.
}\es

\blemma{cap=0} (i) Let $A$ be a \qop subset of $\bdw$ and
$\{u_n\}\sbs\CU(\Gw)$ a \seq of solutions vanishing on $A$. If
$\set{u_n}$ converges then $u=\lim u_n$ vanishes on $A$. In
particular, if $E$ is Borel, $U_E$ vanishes outside $\tl E$.\1
(ii) Let $E$ be a Borel set \sth $\Capq(E)=0$. If $u\in \CU(\Gw)$ and
$u$ vanishes on every $q$-open subset of $E^c=\bdw\sms E$ then
$u\equiv 0$. In particular, $U_E\equiv 0$.\1 (iii) If $\{A_n\}$ is
a \seq of Borel subsets of $\bdw$ \sth $\Capq(A_n)\to 0$ then
$U_{A_n}\to 0$. \es

\bproof (i) Using \rlemma{vanishing} we find that, in proving the
first assertion, we may assume that $\{u_n\}$ is increasing. Now
we can produce an increasing sequence of solutions $\{w_n\}$ such
that, for each $n$, $w_n$ vanishes on some (open) \ngh of $A$ and
$\lim w_n=\lim u_n$. By definition $\lim w_n$ vanishes on $A$.

Let $E$ be a \qcl set. By \rlemma{vanishing}(a) and \rprop{sup
u_n},
 there exists an increasing \seq of solutions $\set{u_n}$
vanishing outside $E$ \sth $U_E=\lim u_n$. Therefore $U_E$
vanishes outside $\tl E$.\1
(ii) Let $A_n$ be open sets \sth $E\sbs A_n$, $ A_n \downarrow$
and $\Capq(A_n)\to 0$. The sets $\tl A_n$ have the same properties
and, by assumption, $u$ vanishes in $(\tl A_n)^c:=\bdw \sms \tl
A_n$. Therefore, for each $n$, there exists a solution $w_n$ which
vanishes on an open \ngh $B_n$ of $(\tl A_n)^c$ \sth $w_n\leq u$
and $w_n\to u$. Hence $w_n\leq U_{K_n}$ where $K_n=B_n^c$ is
compact and $K_n\sbs \tl A_n$. Since  $\Capq(K_n)\to 0$, the
capacitary estimates of \cite{MV4} imply that $\lim U_{K_n}=0$ and
hence $u=0$.\1 (iii)  By definition $U_{A_n}=U_{\tl A_n}$.
Therefore, in view of \rprop{q-top}(iv), it is enough to prove the
assertion when each set $A_n$ is $q$-closed. As before, for each
$n$, there exists
 a solution $w_n$ which vanishes on an open \ngh $B_n$ of $(\tl A_n)^c$
\sth $w_n\leq U_{A_n}$ and $U_{A_n}-w_n\to 0$. Thus $w_n\leq
U_{K_n}$ where $K_n=B_n^c$ is compact and $K_n\sbs \tl A_n$. Since
$\Capq(K_n)\to 0$ it follows that $U_{K_n}\to 0$, which implies
the assertion. \eproof
\blemma{maxsol} Let $E,F$ be Borel subsets of $\bdw$.\1
(i) If $E$, $F$ are  \qcl  then $U_E\wedge U_F= U_{E\cap F}$.\1
(ii) If $E$, $F$ are  \qcl then
\begin{equation}\label{maxsol1}\BAL
& U_E< U_F \iff [\,E\sbsq F \text{ and }\Capq(F\sms E)>0\,],\\
& U_E=U_F \iff E\smq F. \EAL\end{equation} (iii) If $\set{F_n}$ is
a decreasing \seq of \qcl sets then
\begin{equation}\label{maxsol03}
 \lim U_{F_n}= U_F \txt{where} F=\cap F_n.
\end{equation}
(iv) Let $A\sbs \bdw$ be a \qop set and let $u\in\CU(\Gw)$. Suppose
that $u$ vanishes $q$-locally in $A$,
 i.e., for every
point $\gs\in A$ there exists a \qop set $A_\gs$ \sth
$$\gs\in A_\gs\sbs A,\q u\app{A_\gs}0.$$
Then $u$ vanishes on $A$. In particular each solution $u\in \CU(\Gw)$
vanishes on $\bdw\sms \suppq u$. \es \bproof (i) $U_E\wedge U_F$
is the largest solution under $\inf(U_E,U_F)$ and therefore, by
\rdef{maxsol}, it is the largest solution which vanishes outside
$E\cap F$.\1 (ii) Obviously
\begin{equation}\label{maxsol01}
 E\smq F\Lra U_E=U_F,\qq E\sbsq F\iff U_E\leq U_F.
\end{equation}
In addition,
\begin{equation}\label{maxsol02}
  \Capq(F\sms E)>0\Lra U_E\neq U_F.
\end{equation}
Indeed, if $K$ is a compact subset of $F\sms E$ of positive
capacity, then $U_K>0$ and $U_K\leq U_F$ but $U_K\nleq U_E$.
Therefore $U_F=U_E$ implies $F\smq E$.\1 (iii) If $V:= \lim
U_{F_n}$  then $U_F\leq V$. If $U_F<V$ then $\Capq(\suppq V\sms
F)>0.$ But $\suppq V\sbs F_n$ so that $\suppq V\sbs F$ and \consy
$V\leq U_F$.\1
(iv)  First assume that $A$ is a countable union of $\qop$ sets
$\set{A_n}$ \sth $u\app{A_n}0$ for each $n$. Then $u$ vanishes on
$\cup_1^k A_i$ for each $i$. Therefore we may assume that the \seq
$\set{A_n}$ is increasing. Put $F_n=\bdw\sms A_n$. Then $u\leq
U_{F_n}$ and, by (iii),
 $U_{F_n}\downarrow U_F$ where $F=\bdw\sms A$. Thus $u\leq U_F$, i.e., which is equivalent to
 $u\app{A}0$.
 \par We turn to the general case. It is known that the $(\ga,p)$-fine topology
possesses the quasi-Lindel\"{o}f property (see \cite[Sec.
6.5.11]{AH}). Therefore $A$ is covered, up to a set of capacity
zero, by a countable subcover of $\set{A_\gs:\gs\in A}$. Therefore
the previous argument implies that $u\app{A}0$ \eproof
\bth{sup/inf} (a) Let $E$ be a \qcl set. Then,
\begin{equation}\label{sup/inf}\BAL
 U_{ E}&=\inf \{U_D:  E\sbs D\sbs \bdw,\; D\;\text{open}\}\\
 &=\sup \{U_K: K\sbs E,\; K\;\text{compact}\}.
\EAL\end{equation} (b) If $E,F$ are two Borel subsets of $\bdw$
then
\begin{equation}\label{split1}
  U_E= U_{F\cap E}\oplus U_{E\sms F}.
\end{equation}
(c) Let $E,F_n$, $n=1,2,\dots$ be Borel subsets of $\bdw$ and let
$u$ be a positive solution of \req{q-eq}. If either $\Capq(E\Gd
F_n)\to 0$ or $F_n\downarrow E$ then
\begin{equation}\label{Capto0}
  U_{F_n}\to U_E.
\end{equation}
\es \bproof (a) Let $\set{Q_j}$ be a \seq of open sets, decreasing
to a set $E'\smq E$, which satisfies \req{q-closed-app}. Then $\tl
Q_j\downarrow E'$ and, by \rlemma{maxsol}(iii) $U_{Q_j}\downarrow
U_E$.
 This implies
the first equality in \req{sup/inf}. The second equality follows
directly from \rdef{vanishing} (see also \rlemma{vanishing}). \1
(b)  Let $D$, $D'$ be open sets \sth $\widetilde{E\cap F}\sbs D$
and $\widetilde{E\sms F}\sbs D'$ and let $K$ be a compact subset
of $\tl E$. Then
\begin{equation}\label{UKleq}
 U_K\leq U_D +U_{D'}.
\end{equation}
To verify this inequality, let $v$ be a positive solution \sth
$\suppq v\sbs K$ and let $\set{\gb_n}$ be a \seq decreasing to
zero \sth the following limits exist:
$$w=\lim_{n\tin}v_{\gb_n}^D, \qq w'=\lim_{n\tin}v_{\gb_n}^{D'}.$$
(See Notation 1.1 for the definition of $v_\gb^D$.) Then
$$v\leq w+w'\leq U_D +U_{D'}.$$
Since, by \cite{MV4}  $U_K=V_K$, this inequality implies
\req{UKleq}. Further \req{UKleq} and \req{sup/inf}  imply
$$  U_E\leq U_{F\cap E}+ U_{E\sms F}.$$
\par On the other hand, both $U_{F\cap E}$ and $U_{E\sms F}$ vanish outside $\tl E$. \Consy
$U_{F\cap E}\oplus U_{E\sms F}$ vanishes outside $\tl E$ so that
$$  U_E\geq U_{F\cap E}\oplus U_{E\sms F}.$$
This implies \req{split1}.\1 
(c) The previous statement implies,
$$  U_E\leq U_{F_n\cap E}+ U_{E\sms F_n}, \qq U_{F_n}\leq U_{F_n\cap E}+ U_{F_n\sms E}.$$
If $\Capq(E\Gd F_n)\to 0$ , \rlemma{cap=0} implies $\lim U_{E\Gd
F_n}= 0$ which in turn implies  \req{Capto0}.

If $F_n\downarrow E$ then, by \rlemma{maxsol}, $U_{F_n}\to U_E$.
\eproof
\note{Notation 3.2} For any Borel set $E\sbs \bdw$ of positive
$\Capq$-capacity put
\begin{equation}\label{VE}\BAL
 \CV_{mod}(E)&=\{u_\mu: \mu\in \Wqdb,\; \mu(\bdw\sms E)=0\},\1
V_E&=\sup \,\CV_{mod}(E) \EAL\end{equation}

\bth{gencapest} If $E$ is a $q$-closed set, then
\begin{equation}\label{UE=VE1}
  U_E=V_E.
\end{equation}
Thus the maximal solution $U_E$ is $\gs$-moderate. Furthermore
$U_E$ satisfies the capacitary  estimates established in
\cite{MV4} for compact sets, namely:\1 There exist positive
constants $c_1,c_2$ depending only on $q,\, N$ and $\Gw$ \sth, for
every $x\in \Gw $,
  \begin {equation}\label {cap-est}\BAL
c_2\gr(x)&\sum_{m=-\infty}^\infty r_m^{-1-2/(q-1)} C_{2/q,q'}((E\cap S_m(x))/r_m))\leq U_E(x)\leq\\
 c_1\gr(x)&\sum_{m=-\infty}^\infty r_m^{-1-2/(q-1)} C_{2/q,q'}((E\cap S_m(x))/r_m)),
\EAL \end {equation}
 where
$$\gr(x)=\dist(x,\bdw),\q r_m:=2^{-m}, \; S_m(x)=\{y\in \bdw:\,r_{m+1}\leq \abs{x-y}\leq r_{m}\}.$$
Note that, for each point $x$, $S_m(x)=\emptyset$ when
$$\sup_{y\in E}\abs{x-y}< r_{m+1}<r_m<\gr(x).$$
Therefore the sum is finite for each $x\in \Gw$. \es
\Remark
Actually the estimates hold for any Borel set $E$. Indeed, by
definition, $U_E=U_{\tl E}$ and $C_{2/q,q'}((E\cap
S_m(x))/r_m))\sim  C_{2/q,q'}((\tl E\cap S_m(x))/r_m))$. \bproof
Let $\set{E_k}$ be a $q$-stratification of $E$. If $u\in \CV_{mod}(E)$
and $\mu=\tr u$ then $u_\mu=\sup u_{\mu_k}$ where
$\mu_k=\mu\chi\ind{E_k}$. Hence $V_E=\sup V_{E_k}$.
 By \cite{MV4}, $U_{E_k}=V_{E_k}$. These facts and \rth{sup/inf}(c)
imply \req{UE=VE1}. It is known that $U_{E_k}$ satisfies the
capacitary  estimates \req{cap-est}. In addition,
$$C_{2/q,q'}((E_k\cap S_m(x))/r_m))\to C_{2/q,q'}((E\cap S_m(x))/r_m)).$$
Therefore $U_E$ satisfies the capacitary estimates. \eproof

\section{Localization}
\bdef{maxmodsol} {\rm  Let $\mu$ be a positive bounded Borel measure on
$\bdw$ which vanishes on sets of $\Capq$-capacity zero.\1 (a) The
$q$-support of $\mu$ (denoted $\qsupp\mu$) is the intersection of
all \qcl sets $F$ \sth $\mu(\bdw\sms F)=0$.\1 (b) We say that
$\mu$ is concentrated on a Borel set $E$ if $\mu(\bdw\sms E)=0$.
}\es
\blemma{support-mu}
 If $\mu$ is a measure as in the previous definition then,
\begin{equation}\label{sm=su_m}
 \qsupp \mu\smq\,\suppq u_\mu.
\end{equation}
\es \bproof Put $F=\suppq u_\mu$. By \rlemma{maxsol}(iv) $u_\mu$
vanishes on $\bdw\sms F$ and by \rlemma{vanishing} there exists an
increasing \seq of positive solutions $\set{u_n}$ \sth each
function $u_n$ vanishes outside a compact subset of $F$, say
$F_n$, and $u_n\uparrow u_\mu$. If $S_n:=\suppq u_n$ then $S_n\sbs
F_n$ and $\set{S_n}$ increases. Thus $\set{\bar S_n}$ is an
increasing \seq of compact subsets of $F$ and, setting
$\mu_n=\mu\chi\ind{\bar S_n}$, we find $u_n\leq u_{\mu_n}\leq
u_\mu$ so that $u_{\mu_n}\uparrow u_\mu$. This, in turn, implies
(see \cite{MV2})
$$\mu_n\uparrow \mu,\q \qsupp \mu\sbsq \widetilde{\cup_n\bar S_n}\,  \sbs F.$$
If $D$ is a relatively open set and $\mu(D)=0$ it is clear that
$u_\mu$ vanishes on $D$. Therefore $u_{\mu_n}$ vanishes outside
$\bar S_n$, thus outside $\qsupp\mu$.  \Consy $u_\mu$ vanishes
outside $\qsupp\mu$, i.e. $F\sbsq \qsupp \mu$.
\eproof
\bdef{u_D}{\rm
 Let $u$ be a positive solution and $A$ a Borel subset of $\bdw$. Put
\begin{equation}\label{[u]_A}
  [u]_A:=\sup\set{v\in \CU(\Gw): \, v\leq u,\q \suppq v\sbsq \tl A}
\end{equation}
and,
\begin{equation}\label{[u]^A}
  [u]^A:=\sup\set{[u]_F: F\sbsq A,\;F \text{ \qcl}}.
\end{equation}
Thus $[u]_A=u\wedge U_A$, i.e., $[u]_A$ is the largest solution
under  $\inf (u,U_A)$.
\par Recall that, if $A$ is \qop and $u\in C(\bdw)$,
$u^A_\gb$ denotes the solution of \req{q-eq} in $\Gw'_\gb$ which equals
$u\chi\ind{\Gs_\gb(A)}$ on $\Gs_\gb$.
\par If $\lim_{\gb\to 0} u_\gb^A$ exists the limit will be denoted by $u^A$.
} \es
\bth{lim[u]} Let $u\in \CU(\Gw)$.\1 (i) If $E\sbs\bdw$ is \qcl then,
\begin{equation}\label{sup/inf1}
 [u]_E=\inf \{[u]_D:  E\sbs D\sbs \bdw,\; D\;\text{open}\}.
\end{equation}
(ii) If $E,F$ are two Borel subsets of $\bdw$  then
\begin{equation}\label{split2}
[u]_E\leq [u]_{ F\cap  E}+ [u]_{E\sms F}
\end{equation}
and
\begin{equation}\label{double[}
  [[u]_E]_F=[[u]_F]_E=[u]_{ E\cap F}.
\end{equation}
(iii) Let $E,F_n$, $n=1,2,\dots$ be Borel subsets of $\bdw$.  If
either $\Capq(E\Gd F_n)\to 0$ or $F_n\downarrow E$ then
\begin{equation}\label{Capto0'}
 [u]_{F_n}\to [u]_E.
\end{equation}
\es \bproof (i) Let $\CD=\set{D}$ be the family of sets in
\req{sup/inf1}. By \req{sup/inf}  (\wrto the family $\CD$ )
\begin{equation}\label{sup/inf2}
  \inf(u,U_E)=\inf(u,\inf\ind{D\in\CD}U_D)=\inf\ind{D\in\CD}\inf(u,U_D)\geq \inf\ind{D\in\CD}[u]_D.
\end{equation}
Obviously $$[u]_{D_1}\wedge [u]_{D_2}\geq [u]_{D_1\cap D_2}.$$
 (In fact we have equality but that is not needed here.)
Therefore, by \rprop{sup u_n}, the function
$v:=\inf\ind{D\in\CD}[u]_D$ is a solution of \req{q-eq}. Hence
\req{sup/inf2} implies
 $[u]_E\geq v$. The opposite inequality is obvious.\1
(ii)  If $E$ is compact \req{split2} is proved in the same way as
\rth{sup/inf}(b). In general, if $\set{E_n}$ is a $q$
stratification of $\tl E$,
$$[u]_{E_n}\leq [u]_{F\cap E_n}+ [u]_{E_n\sms F}\leq [u]_{F\cap E}+ [u]_{E\sms F}.$$
This inequality and \rth{sup/inf}(c) imply \req{split2}.
\par Put $A=\tl E$ and $B=\tl F$. It follows directly from the definition that,
$$[[u]_A]_B\leq inf(u,U_A,U_B).$$
The largest solution dominated by $u$ and vanishing on $A^c\cup
B^c)$ is $[u]_{A\cap B}$. Thus
$$[[u]_A]_B\leq [u]_{A\cap B}.$$
On the other hand
$$[u]_{A\cap B}=[[u]_{A\cap B}]_B\leq [[u]_A]_B.$$
This proves \req{double[}.
(iii) By \req{sup/inf2}
$$  [u]_E\leq [u]_{F_n\cap E}+ [u]_{E\sms F_n}, \qq [u]_{F_n}\leq [u]_{F_n\cap E}+ [u]_{F_n\sms E}.$$
If $\Capq(E\Gd F_n)\to 0$ , \rlemma{cap=0} implies $\lim [u]_{E\Gd
F_n}= 0$ which in turn implies \req{Capto0'}.

\par If $F_n\downarrow E$ then, by \rlemma{maxsol}, $U_{F_n}\to U_E$. If $u$ is a positive solution then
$$\inf(u,U_E)=\inf(u,\inf_n U_{F_n})=\inf_n\inf(u,U_{F_n})\geq \inf_n[u]_{F_n}.$$
Since $\set{F_n}$ decreases $w=\inf_n[u]_{F_n}$ is a solution.
Hence $[u]_E\geq w$. The opposite inequality is obvious; hence
$[u]_E=\lim [u]_{F_n}$. \eproof

\blemma{convergence2} Let $u$ be a positive solution of \req{q-eq}
and put $E=\suppq u$.\1
(i) If $D$ is a $q$-open set \sth $E\sbsq D$ then
\begin{equation}\label{conv1}
[u]^D= \lim_{\gb\to 0}u_\gb^D=[u]_D=u.
\end{equation}
(ii) If $A$ is a \qop subset of $\bdw$,
\begin{equation}\label{conv1'}
  u\app{A}0 \iff u^Q=\lim_{\gb\to 0}u_\gb^{Q}=0 \forevery Q \; \text{\rm  \qop}: \tl Q\sbsq A.
\end{equation}
(iii) Finally,
\begin{equation}\label{conv1*}
 u\app{A}0 \Llra [u]^{A}=0
\end{equation}
\es \bproof \note{Case 1: $E$ is  closed}  Since $u$ vanishes in $A:=\bdw\sms E$,
it follows that $u\in C(\Gw\cup A)$ and $u=0$ on $A$. If, in
addition, $D\sbs \bdw$ is an {\em open} \ngh of $E$ then
$$\int_{\Gs_\gb(D^c)}udS\to 0$$
so that
\begin{equation}\label{conv1-a}
 \lim_{\gb\to 0}u_\gb^{D^c}=0.
\end{equation}
Since
$$u_\gb^{D}\leq u\leq u_\gb^{D}+ u_\gb^{D^c} \txt{in} \Gw'_\gb$$
it follows that
\begin{equation}\label{conv1-b}
 u=\lim u_\gb^D.
\end{equation}

\par If we assume only that $D$ is $q$-open and $E\sbsq D$ then, for every $\ge>0$, there exists an
open set $O_\ge$ \sth $D\sbs O_\ge$, $E\sbs O_\ge$ and $\Capq(O'_\ge)<\ge$ where $O'_\ge=O_\ge\sms
D$. It follows that
$$u_\gb^{ O_\ge}-u_\gb^D\leq U_{\Gs_\gb(O'_\ge)} \txt{in} \Gw'_\gb$$
and $\lim_{\ge\to 0} U_{\Gs_\gb(O'_\ge)}=0$ uniformly \wrto $\gb$.
 Since $\lim_{\gb\to 0}u_\gb^{O_\ge}=u$ it follows that \req{conv1-b} holds.
The same argument shows that \req{conv1-a} remains valid.

Now \req{conv1-b} implies \req{conv1}. Indeed
  $$u=\lim u_\gb^D\leq [u]_D\leq u.$$
  Hence $u=[u]_D$. If $Q$ is a \qop set \sth
 $E\sbsq Q\sbs \tl Q\sbsq D$ then $u=[u]_Q\leq [u]^D$. Hence $u=[u]^D$.

 \par In addition \req{conv1-a} implies \req{conv1'} in the direction $\Lra$. Assertion
\req{conv1'} in the opposite direction is a \cons of
\rlemma{closedngh} and \rlemma{maxsol} (iv).
\ \1 \note{Case 2} We consider the general case when $E$ is
$q$-closed. Let $\{E_n\}$ be a stratification of $E$ so that
$\Capq(E\sms E_n)\to 0$. If $D$ is  $q$-open and $E\sbsq D$ then, by
the first part of the proof,
\begin{equation}\label{sup/inf3}
  \lim_{\gb\to 0}([u]_{E_n})_\gb^D=[u]_{E_n}.
\end{equation}
By \req{split2}
\begin{equation}\label{sup/inf4}
u_\gb^D\leq ([u]_{E_n})_\gb^D+ ([u]_{E\sms E_n})_\gb^D.
\end{equation}
Let $\set{\gb_k}$ be a \seq decreasing to zero \sth the following
limits exist
$$w:=\lim_{k\to 0}u_{\gb_k}^D, \q w_n:=\lim_{k\to 0}([u]_{E\sms E_n})_{\gb_k}^D,\q n=1,2,\cdots\,.$$
Then, by \req{sup/inf3} and \req{sup/inf4},
$$[u]_{E_n}\leq w\leq [u]_{E_n}+ w_n\leq [u]_{E_n}+ U_{E\sms E_n}.$$
Further, by \req{Capto0},
$$[u]_{E_n}\to [u]_E=u,\q U_{E\sms E_n}\to 0.$$
Hence $w=u$. This implies \req{conv1-b} which in turn implies
\req{conv1}.
\par To verify \req{conv1'} in the direction $\Lra$ we apply \req{sup/inf4} with $D$ replaced by $Q$.
 We obtain,
$$u_\gb^{Q}\leq ([u]_{E_n})_\gb^{Q}+ ([u]_{E\sms E_n})_\gb^{Q}.$$
By the first part of the proof
$$\lim_{\gb\to 0}([u]_{E_n})_\gb^{Q}=0.$$
Let $\set{\gb_k}$ be a \seq decreasing to zero \sth the following
limits exist:
$$\lim_{k\tin}([u]_{E\sms E_n})_{\gb_k}^{Q},\q n=1,2,\ldots,\q\lim_{k\tin}u_{\gb_k}^{Q}.$$
Then
$$\lim_{k\tin}u_{\gb_k}^{Q}\leq \lim_{k\tin}([u]_{E\sms E_n})_{\gb_k}^{Q} \leq U_{E\sms E_n}.$$
Since $U_{E\sms E_n}\to 0$ we obtain \req{conv1'} in the direction
$\Lra$. The assertion in the opposite direction is proved as in
Case 1. This completes the proof of (i) and (ii).
\par Finally we prove (iii). First assume that $u\app{A}0$.
If $F$ is a \qcl set \sth $F\sbsq A$ then there exists a \qop set
$Q$ \sth $F\sbsq \tl Q\sbsq A$. Therefore, applying \req{conv1} to
$v:=[u]_F$ and  using \req{conv1'} we obtain
$$ v=\lim v_\gb^{Q}\leq \lim u_\gb^{Q}=0.$$
In view of \rdef{u_D} this implies that $[u]^A=0$.
\par Secondly assume that  $[u]^A=0$. Then $[u]_Q=0$ whenever $\tl Q\sbsq A$. If $Q$ is a \qop set
\sth $\tl Q\sbsq A$ then $[u]_Q=0$ and hence $u\app{Q}0$. Applying
once again \rlemma{closedngh} and \rlemma{maxsol} (iv) we conclude that
$u\app{A}0$. \eproof
\bdef{loctr1}{\rm
 Let $u,v$ be positive solutions of \req{q-eq} in $\Gw$  and let $A$ be a \qop subset of $\bdw$.
  We say that $u=v$ on $A$ if $u\ominus v$ and $v\ominus u$ vanish on $A$ (see Notation 3.1). This relation is denoted by
$u\underset{A}{\approx}v$. }\es
\bth{uappv} Let $u,v\in \CU(\Gw)$  and let $A$ be a \qop subset of $\bdw$. Then,
\begin{equation}\label{uappv1}
 u\app{A}v \Llra \lim_{\gb\to 0} \abs{u-v}_\gb^Q=0,
\end{equation}
for every \qop set $Q$ \sth $\tl Q\sbsq A$ and
\begin{equation}\label{uappv2}
 u\app{A}v \iff [u]_F=[v]_F,
\end{equation}
for every \qcl set $F$ \sth $F\sbsq A$.
\es
\bproof By definition,
$u\app{A}v$ is equivalent to $u\ominus v\app{A}0$ and $v\ominus
u\app{A}0$. Hence, by \rlemma{convergence2} (specifically
\req{conv1*}),
\begin{equation}\label{u-v}
 [u\ominus v]_F=0,\q [v\ominus u]_F=0,
\end{equation}
for every \qcl set $F\sbsq A$. Therefore, if $\tl Q\sbsq A$,
\rlemma{convergence2} implies that
$$((u-v)_+)_\gb^Q\to 0, \q ((v-u)_+)_\gb^Q\to 0.$$
(Recall that $u\ominus v$ is the smallest solution which dominates
the subsolution $(u-v)_+$.)
 This implies \req{uappv1} in the direction $\Lra$;
the opposite direction is a \cons of \rlemma{cap=0}.
\par We turn to the proof of \req{uappv2}. For any two positive solutions $u,v$ we have
\begin{equation}\label{u-v,v-u}
 u+(v-u)_+\leq v+(u-v)_+\leq v+u\ominus v.
\end{equation}
\bcom
which  implies
\begin{equation}\label{uominusv}
u+v\ominus u\leq v+u\ominus v.
\end{equation}
To verify this implication, let $w_{1,\gb}$ (resp. $w_{2,\gb}$) be the solution of \rqq in $\Gw'_\gb$
with boundary data $(v-u)_+$ (resp. $(u-v)_+$) on $\Gs_\gb$. Since $(v-u)_+$ and $(u-v)_+$)
are subsolutions, $w_{1,\gb}\downarrow v\ominus u$ and $w_{2,\gb}\downarrow u\ominus v$.
Obviously,
$$u_\gb+ w_{1,\gb}\leq v_\gb+w_{2,\gb},$$
where $u_\gb=u\chr{\Gw'_\gb}$.
\end{comment}
If $F$ is a \qcl set and $Q$ a \qop set \sth $F\sbsq Q$ then,
\begin{equation}\label{uFv}
  [u]_F\leq [v]_Q+ [u\ominus v]_Q.
\end{equation}
To verify this inequality we observe that, by \req{u-v,v-u},
$$[u]_F\leq [v]_Q+[v]_{Q^c}+[u\ominus v]_Q + [u\ominus v]_{Q^c}.$$
The subsolution $w:=([u]_F -([v]_Q+ [u\ominus v]_Q))_+$ is dominated
by the  supersolution $[v]_{Q^c} + [u\ominus v]_{Q^c}$ which vanishes on $Q$. Therefore $w$ vanishes on
$Q$. Since the boundary support of $[w]_\dag$ is contained in $F$ it follows that $[w]_\dag\equiv 0$
so that $w\equiv0$.

 \par If $u\app{A}v$ and $F\sbsq Q\sbs \tl Q\sbsq A$ then \req{uFv} and \req{u-v} imply,
 $$[u]_F\leq [v]_Q.$$ Choosing a decreasing \seq of \qop sets $\set{Q_n}$ \sth $\cap Q_n\smq F$
 we obtain $[u]_F\leq \lim [v]_{Q_n}=[v]_F$.
Similarly, $ [v]_F\leq [u]_F$ and hence equality.

\par Next assume that $ [v]_F=[u]_F$ for every \qcl set $F\sbsq A$. If $Q$ is a \qop set \sth
$\tl Q\sbsq A$ we have,
$$u\ominus v\leq ([u]_Q\oplus[u]_{Q^c})\ominus [v]_Q\leq [u]_{Q^c},$$
because $[u]_Q=[v]_Q$. This implies that $u\ominus v$ vanishes on $Q$. Since this holds for every $Q$ as above it
follows that $u\ominus v$ vanishes on $A$. Similarly $v\ominus u$ vanishes on $A$.
\eproof
\bcor{uappv} If $A$ is a \qop subset of $\bdw$, the relation
$\app{A}$ is an equivalence relation in $\CU(\Gw)$.
\es

\bproof This is
an immediate \cons of \req{uappv1}.
\eproof

\section{The precise boundary trace}
\subsection{The regular boundary set} We define the regular boundary set of a positive solution of \rqq
and present some conditions for the regularity of a \qop set.
\bdef{trace}{\rm
 Let $u$ be a positive solution of \rqq.
\begin{description}
  \item[{\rm a}] Let $D\sbs \bdw$ be a \qop set \sth $\Capq(D)>0$. $D$ is {\em pre-regular} \wrto $u$ if 
\begin{equation}\label{regularity}
 \int_{\Gw}[u]_F^q\gr dx<\infty \forevery F\sbsq D,\;F \text{ \qcl}.
\end{equation}
  \item[{\rm b}] An arbitrary Borel set $E$ is {\em regular} if there exists a pre-regular  set $D$ \sth
  $\tl E\sbsq D$.
  \item[{\rm c}] A set $D\sbs \bdw$ is $\gs$-{\em  regular} if it is the union of a countable family
  of pre-regular sets.
 \item[{\rm d}]  The union of all \qop regular sets is called {\em the  regular boundary set of} $u$,
and is denoted by $\CR(u)$. The set $\CS(u)=\bdw\sms \CR(u)$ is called
{\em the singular boundary set} of $u$. A point $P\in \CR(u)$ is
called a regular boundary point of $u$; a point $P\in \CS(u)$ is
called a singular boundary point of $u$.
 \end{description}
 }\es
 \Remark The property of regularity of a set is preserved under the equivalence relation $\smq$.
However note that {\em a point is regular \ifif it has a \qop regular \ngh\/}.
\blemma{pre-reg}
 If $D$ is a \qop pre-regular set then every point $\gx\in D$ is a regular point. Furthermore there
exists a \qop regular set $Q$ \sth
\begin{equation}\label{reg/sing}
\gx\in Q\sbs \tl Q\sbs D.
\end{equation}
If $F$ is a regular \qcl set  then there exists a regular \qop set $Q$ \sth $F\sbsq Q$.
\es
\bproof  By \rlemma{closedngh}, for every $\gx\in D$, there
exists a \qop set $Q$ \sth \req{reg/sing} holds. Therefore $Q$ is
a regular set and $\gx$ is a regular point. The last assertion is a consequence of \rlemma{e-ngh}.
\eproof
 \bdef{decomposition} {\rm
 Let $u$ be a positive solution of \req{q-eq} and let $\set{Q_n}$ be an increasing \seq
  of regular \qop sets. If $\tl Q_n\sbsq Q_{n+1}$ we say that $\set{Q_n}$ is
  {\em  a  regular \seq} relative to $u$.

  If $Q$ is a \qop set, $\set{Q_n}$ is
   a  regular \seq relative to $u$, $Q_n\sbs Q$ and $ Q_0:=\cup_{n=1}^\infty Q_n\smq Q$
  we refer to $Q_0$  as a {\em proper representation} of  $Q$ and to $\set{Q_n}$ as a {\em regular decomposition}
  of $Q$, relative to $u$.
}\es
\blemma{decomposition}  Let $u\in \CU(\Gw)$.
 A \qop set $Q\sbs \bdw$ is \gsreg \ifif it has a
proper representation relative to $u$. In particular every
pre-regular set has a proper representation.
\es
\bproof The 'if' direction
follows immediately from the definition. Now suppose that $Q$ is
\gsreg. Then $Q=\cup_1^\infty E_n$ where $E_n$ is  \qop and
pre-regular, $n=1,2,\cdots\,$.
  By \rlemma{set-approx},
 each set $E_n$ can be represented
(up to a set of capacity zero) as a countable union of  \qop sets
$\{A_{n,j}\}_{j=1}^\infty$ \sth $\tl A_{n,j}\sbsq A_{n,j+1}\sbsq E_n$.
We may assume that $A_{n,j}\sbs E_n$; otherwise we replace it by $A_{n,j}\cap E_n$.
Put
$$Q_n=\cup_{k+j=n} A_{k,j}.$$
If $k+j=n$   then $\tl A_{k,j}\sbsq A_{k,j+1}\sbs Q_{n+1}$. Hence
$$\tl Q_n\sbsq Q_{n+1},\q Q_0:=\cup Q_n\smq Q.$$
\eproof
\bth{maintrace} Let $D$ be a \qop set \sth $\Capq(D)>0$.\1
\num{i} Suppose
that
\begin{equation}\label{regcondition}
  \liminf_{\gb\to 0} \int_{\Gw'_\gb} (u_\gb^D)^q(\gr-\gb)\, dx<\infty.
\end{equation}
Then $D$ is pre-regular.\1
(ii) Suppose that $D$ is a pre-regular set. Then there exists a
Borel measure $\mu$ on $D$ \sth, for every  \qcl set  $E\sbsq D$,
\begin{equation}\label{regcond4}
  \tr[u]_E= \mu\chr{E}.
\end{equation}
\es
\bcom Furthermore,
\begin{equation}\label{regcondition'}
 u^D:=\lim_{\gb\to 0} u^D_\gb
\end{equation}
is well defined and
\begin{equation}\label{regcondition''}
 [u^D]_E=[u]_E
\end{equation}
for every \qcl set $E\sbsq D$.
\end{comment}
\bproof (i) Let $\{\gb_n\}$ be a \seq decreasing to zero \sth
\begin{equation}\label{regcond1}
   \int_{\Gw'_{\gb_n}} (u_{\gb_n}^D)^q(\gr-\gb_n)\, dx\leq C <\infty, \forevery n.
\end{equation}
By extracting a \sseq if necessary we may assume that
$\{u_{\gb_n}^D\}$ converges locally uniformly in $\Gw$ to a
solution $w$. Then, by \rlemma{convergence2}, if  $E$ is  \qcl and
$E\sbsq D$,
\begin{equation}\label{regcond2}
 [u]_E=\lim_{\gb\to 0} ([u]_E)_\gb^D\leq \lim_{n\tin} u_{\gb_n}^D=: w.
\end{equation}
 By \req{regcond1} and Fatou's lemma,
$$ \int_{\Gw} w^q\gr\, dx\leq C <\infty.$$
Hence, by \req{regcond2},
\begin{equation}\label{regcond3}
  \int_{\Gw} ([u]_E)^q\gr\, dx\leq C <\infty.
\end{equation}
Thus $D$ is pre-regular.\1
(ii) By \rlemma{decomposition}, $D$
possesses a regular decomposition $\set{D_j}$. Put
$$w_j=[u]_{D_j}.$$
Then $\set{w_j}$ is increasing and its limit is a solution
$w_0\leq w$ with $w$ as defined in \req{regcond2}. Thus $w_0$ is a
moderate solution.
 If $E\sbsq D$ is a \qcl set  then, by \req{double[},
$$[u]_{E\cap \tl D_j}=[w_j]_E\leq [w_0]_E\leq [u]_E.$$
By \rlemma{q-compact}, for every $k\in \BBN$ there exists an open
set $O_k$ and a natural number $j_k$  \sth $\Capq(O_k)<1/k$ and
$E\sms O_k\sbsq D_{j_k}$. By \rth{sup/inf}
$$[u]_E\leq [u]_{E\cap \tl D_{j_k}}+[u]_{O_k}.$$
Since $[u]_{O_k}\to 0$ we conclude that
\begin{equation}\label{w0E}
[w_0]_{E}= [u]_{E} \forevery E\sbsq D:\; E \text{ \qcl.}
\end{equation}
 If $w_0$ is moderate then
$\tr [w_0]_{E}=\mu\chr{E}\tr w_0$, which implies \req{regcond4}.

We turn to the case where $w_0$ is not moderate. The solution $w_j$ is moderate and we denote
$$\mu_j=\tr w_j,\q \mu=\lim\mu_j.$$
 By \req{double[}, $w_j=[w_{j+k}]_{D_j}$. Therefore $\mu_j=\mu_{j+k}\chr{D_j}=\mu\chr{D_j}$.
 Therefore if $E$ is \qcl and $E\sbsq D_j$ for some $j$, \req{regcond4} holds with $\mu$ as defined above.
 If $E$ is \qcl and $E\sbsq D$ then $E\smq E':=\cup (E\cap D_j)=\cup (E\cap \tl D_j)$.
 Put $E_j:=E\cap \tl D_j$. It follows that
 $$\tr [u]_{E_j}=\mu\chr{E_j}\uar \mu\chr{E'}.$$
 Since $D$ is pre-regular,  $[u]_E$ is moderate. Put $E'_j=E\sms D_j$ and observe that
 $\cap_1^\infty E'_j$ is a set of capacity zero so that (by \rlemma{maxsol}) $U_{E'_J}\dar 0$ and hence
 $\lim [u]_{E'_j}\dar 0$. Since
 $$[u]_E\leq [u]_{E_j}+[u]_{E'_j} \txt{and} [u]_{E'_j}\dar 0$$
 we conclude that
 $$\tr [u]_E\leq \lim \tr [u]_{E_j}=\mu\chr{E'}.$$
On the other hand
$$\tr[u]_E\geq \tr[u]_{E_j}\to\mu\chr{E'}=\mu\chr{E}.$$
 This implies \req{regcond4}.
\eproof
\bcor{maintrace} Let $D$ be a \qop set \sth $\Capq(D)>0$. Suppose
that, for every \qop set $Q$ \sth $\tl Q\sbsq D$,
\begin{equation}\label{loc-regcond}
  \liminf_{\gb\to 0} \int_{\Gw'_\gb} (u_\gb^Q)^q(\gr-\gb)\, dx<\infty.
\end{equation}
Then $D$ is pre-regular.
\es
\bproof This is an immediate \cons of
\rlemma{e-ngh} and \rth{maintrace}.
\eproof
\bcom \blemma{CR(u)} If $u$ is a positive solution of \req{q-eq}
then $\CR(u)$ is \gsreg and \consy possesses a proper
representation.
\par If $\CR_0(u)$ is a proper representation then
$\CR_0(u)\sbs\CR(u)$. \es
\bproof  By \cite[Sec. 6.5.11]{AH} the $(\ga,p)$-fine topology
possesses the quasi-Lindel\"{o}f property.  This  implies that
$\CR(u)$  is \gsreg. By \rlemma{decomposition} $\CR(u)$ has a
proper representation $\CR_0(u)=\cup Q_n$. Since each set $Q_n$ is
regular, $\CR_0(u)\sbs \CR(u)$.
\end{proof}
\end{comment}
\subsection{Behavior of the solution at the boundary}
In this subsection we provide a characterization of regular and singular  boundary points
 of a positive solution by the limiting behavior of the solution as it approaches the boundary in a \qop
 \ngh of each point.
  \bth{loc-tr2} Let $u$ be a positive solution of \req{q-eq} in $\Gw$ and let $Q$ be a \qop subset
  of $\bdw$ of positive capacity.\1
  (i) If  $\gx\in \CS(u)$ then,  for every nontrivial \qop \ngh $Q$ of $\gx$,
\begin{equation}\label{limsing}
 \int_{\Gs_\gb(Q)}u dS\tin, \txt{as} \gb\to 0.
\end{equation}
(ii)  If $\gx\in \CR(u)$, there exists a \qop regular set $D$ \sth $\gx\in D$. Further there exists
 a \qop set $Q$ \sth $\gx\in Q\sbs \tl Q\sbs D$.
 \Consy
\begin{equation}\label{limreg}
\sup_{0<\gb<\gb_0} \int_{\Gs_\gb(Q)}u dS<\infty,\qq  u^Q=\lim_{\gb\to0}u_\gb^Q \text{ exists}
\end{equation}
and $u^Q$ is moderate.
If $\mu\ind{Q}:=\tr [u]_Q$ then, for every  \qcl set  $E\sbsq Q$,
\begin{equation}\label{regcond-bis}
  \tr[u]_E= \mu\ind{Q}\chr{E}.
\end{equation}
\bcom
 $u^Q_\gb\modcon \mu\ind{\CR}\chr{Q}$ as $\gb\to 0$.
In other words, $u^Q=\lim u_\gb^Q$ exists and $u^Q$ is a moderate
solution \sth
\begin{equation}\label{L1bound}
 \tr u^Q= \mu\ind{\CR}\chr{Q}.
\end{equation}
\end{comment}
\es
\bproof (i) If $Q$ is a \qop set for which \req{limsing} does not
hold, there exists a \seq $\{\gb_n\}$ converging to zero \sth
$$ \int_{\Gs_{\gb_n}(Q)}u dS\to \ga<\infty.$$
This implies that there exists a constant $C$ \sth \req{regcond1}
holds (with $D$ replaced by $Q$). By \rth{maintrace}, $Q$ is
pre-regular and by \rlemma{pre-reg} every point in $Q$ is a
regular point. Therefore, if $\gx\in \CS(u)$, \req{limsing} holds.\1
(ii) By \rdef{trace} there exists a \qop regular \ngh $D$ of $\gx$.
By \rlemma{closedngh} there exists a  \qop set $Q$ \sth $\gx\in Q\sbs \tl Q\sbs D$.
By \rth{lim[u]} and \rlemma{convergence2},
$$u\leq [u]_D+[u]_{D^c} \txt{and} \lim_{\gb \to0}\int_{\tl Q} [u]_{D^c}(\gb,\cdot) dS=0.$$
 Therefore
 $$\limsup_{\gb\to 0}\int_Q u(\gb,\cdot)dS\leq \limsup_{\gb\to 0}\int_Q [u]_D(\gb,\cdot)dS<\infty$$
 so that \req{regcondition} holds.
In view of this fact, \rth{maintrace} and the arguments in its proof imply assertion (ii).
\bcom
(ii) By \rlemma{pre-reg} and \rth{loc-tr1}, if $\gx\in\CR(u)$,
there exists   a \qop regular set $D$ \sth $\gx\in D\sbs \tl
D\sbsq \CR(u)$ and $\mu\ind{\CR}(\tl D)<\infty$.
 The solution
$v=[u]_D$ is moderate.  By \rth{loc-tr1},
$u\app{\CR(u)}v\ind{\CR}$ and by \req{consist} $\tr v=\chr{\tl
D}\mu\ind{\CR}$. Hence $v\app{D}v\ind{\CR}$ and \consy, by
\rcor{uappv}, $u\app{D}v$. Therefore, by \rth{uappv},
\begin{equation}\label{loctr1}
\lim_{\gb\to 0} \abs{u-v}_\gb^Q=0,
\end{equation}
for every \qop set $Q$ \sth $\tl Q\sbsq D$.
By \rprop{mod-tr}
$$\lim_{\gb\to 0}\int_{\Gs_\gb}v^\gb_*dS=\mu\ind{\CR}(\tl D)<\infty.$$
This fact and \req{loctr1} imply \req{L1bound}.
\end{comment}
\eproof

\subsection{$q$-perfect measures}
\bdef{cap-reg} {\rm
Let $\mu$ be a positive Borel measure,
not necessarily bounded, on $\bdw$.\1
(i) We say that $\mu$ is {\em essentially  absolutely continuous} relative to  $\Capq$ if the
following condition holds:
\par If $Q$ is a \qop set and $A$ is a Borel set \sth $\Capq(A)=0$ then
 $$\mu(Q\sms A)=\mu(Q).$$
This relation will be denoted by $\mu\ppf\Capq$.\1
(ii) $\,\mu$ is {\em regular relative to $q$-topology} if, for every Borel set $E\sbs \bdw$,
\begin{equation}\label{cap-reg}\BAL
 \mu(E)&=\inf\set{\mu(D): E\sbs D\sbs \bdw,\; D \text{ \qop}}\\
&= \sup\set{\mu(K): K\sbs E,\; K \text{ compact}}.
\EAL\end{equation}
\par $\,\mu$ is {\em outer regular relative to $q$-topology} if the first equality in \req{cap-reg} holds.\1
(iii)  A positive Borel measure is called $q$-{\em perfect} if it is essentially  absolutely
continuous relative to $\Capq$ and outer regular relative to $q$-topology.
The space of $q$-perfect Borel measures is denoted by $\BBM_q(\bdw)$.
}\es
\blemma{cap-reg} If $\mu\in\BBM_q(\bdw)$ and
$A\sbs \bdw$ is a non-empty Borel set \sth $\Capq(A)=0$ then
\begin{equation}\label{negligibleset}
\mu(A)=
  \begin{cases}
    \infty & \text{if $\mu(Q\sms A)=\infty\forevery Q$ \qop \ngh of $A$}, \\
    0 & \text{otherwise}\,.
  \end{cases}
\end{equation}
\par If $\mu_0$ is an essentially \abc positive Borel measure on $\bdw$ and $Q$ is a \qop
set \sth $\mu_0(Q)<\infty$ then $\mu_0\rest{Q}$ is \abc \wrto $\Capq$ in the strong sense, i.e.,
if $\set{A_n}$ is a \seq of Borel subsets of $\bdw$,
$$\Capq(A_n)\to 0 \Lra \mu_0(Q\cap A_n)\to 0.$$
\par Let $\mu_0$ be an essentially \abc positive Borel measure on $\bdw$.  Put
\begin{equation}\label{regext}
\mu(E):=\inf\set{\mu_0(D): E\sbs D\sbs \bdw,\; D \text{ \qop}},
\end{equation}
for every Borel set $E\sbs \bdw$. Then
\begin{equation}\label{capreg0}\BAL
(a)\q& \mu_0\leq \mu, \qq \mu_0(Q)=\mu(Q) \forevery Q \q \text{\qop}\\
(b)\q& \mu\rest{Q}=\mu_0\rest{Q}  \txt{for every \qop set $Q$  \sth $\mu_0(Q)<\infty$}.
 \EAL
\end{equation}
Finally $\mu$ is $q$-perfect; thus $\mu$ is the smallest measure in $\BBM_q$ which dominates $\mu_0$.
\es

\bproof
The first assertion follows immediately from the definition of $\BBM_q$.
We turn to the second assertion.
If $\mu_0$ is an essentially \abc positive Borel measure on $\bdw$ and $Q$ is a \qop
set \sth $\mu_0(Q)<\infty$ then $\mu_0\chr{Q}$ is a bounded Borel measure which vanishes on sets
of $\Capq$-capacity zero. If $\set{A_n}$ is a \seq of Borel sets \sth $\Capq(A_n)\to 0$
and $\mu_n:=\mu_0\chr{Q\cap A_n}$ then $u_{\mu_n}\to 0$ locally uniformly and
$\mu_n\rightharpoonup 0$ weakly \wrto $C(\bdw)$. Hence $\mu_0(Q\cap A_n)\to 0$. Thus $\mu_0$ is
\abc in the strong sense relative to $\Capq$.
\par Assertion \req{capreg0} (a)
follows from \req{regext}. It is also clear that $\mu$, as defined by \req{regext}, is a measure.
Now if $Q$ is a \qop set \sth
$\mu_0(Q)<\infty$ then $\mu(Q)<\infty$ and both $\mu_0\rest{Q}$ and $\mu\rest{Q}$ are regular
relative to the induced Euclidean topology on $\bdw$.
Since they agree on open sets, the regularity implies \req{capreg0} (b).
\par If $A$ is a Borel set \sth $\Capq(A)=0$ and $Q$ is a \qop set then $Q\sms A$ is \qop and \consy
$$\mu(Q)=\mu_0(Q)=\mu_0(Q\sms A)=\mu(Q\sms A).$$
Thus $\mu$ is essentially \abc. It is obvious by its definition that $\mu$ is outer regular \wrto
$\Capq$. Thus $\mu\in\BBM_q(\bdw)$.
\eproof
\subsection{The boundary trace on the regular set}
\par First  we describe some properties  of moderate
solutions. In this connection it is convenient to
introduce a related term: a solution is {\em strictly
moderate} if
\begin{equation}\label{stricthar}
  \abs{u}\leq v, \q v\text{ harmonic},\q \int_\Gw v^q\gr\,dx<\infty.
\end{equation}
A positive solution $u$ is strictly moderate \ifif $\tr u\in \Wqdb$ (see \cite{MV3}).

\note{Notation 5.1} Let $\Gp:\Gw_{\gb_0}\to\bdw$ be the mapping given by $\Gp(x)=\gs(x)$ (see Notation 1.1) and put
$\Gp_\gb:=\Gp\big|\ind{\Gs_\gb}$.

\begin{description}
  \item[1]If $\gf$ is a function defined on $\bdw$ put  $\gf^*:=\gf\circ\Gp$. This function is
   called the {\em normal lifting} of  $\gf$ to $\Gw_{\gb_0}$.
   Similarly, if $\gf$ is defined on a set $Q\sbs \bdw$, $\gf^*$ is the normal lifting of
  $\gf$ to $\Gw_{\gb_0}(Q)$.
  \item[2] If $\vgf$ is a function defined on $\Gs_\gb$ we define the {\em normal projection} of
   $\vgf$ onto $\bdw$ by
  $$\vgf_*^\gb(\gx)=\vgf(\Gp_\gb^{-1}(\gx)), \forevery \gx\in\bdw, \q \Gp_\gb=\Gp\big|\ind{\Gs_\gb}:\Gs_\gb\mapsto\bdw.$$
If $v$ is a function defined on $\Gw_{\gb_0}$ then $v_*^\gb$
denotes the normal projection of $v(\gb,\cdot)$ onto $\bdw$, for $\gb\in
(0,\gb_0)$.
 \end{description}

\bprop{mod-tr} Let $u$ be a  moderate solution of \rqq, not
necessarily positive. Then:\2 (i) $u\in L^1(\Gw)\cap L^q(\Gw;\gr)$
and $u$ possesses a boundary trace $\tr u$ given by a bounded
Borel measure $\mu$ which is attained in the sense of weak
convergence of measures:
\begin{equation}\label{limtr}
  \lim_{\gb\to 0} \int_{\Gs_\gb}u \gf^* dS=\lim_{\gb\to 0} \int_{\bdw}u_*^\gb \gf dS= \int_{\bdw}\gf\,d\mu,
\end{equation}
for every $\gf\in  C(\bdw)$.\1 (ii) A  {\em bounded}  Borel
measure $\mu$ is the boundary trace of a
 solution of \rqq \ifif \ it is \abc relative to $\Capq$. When this is the case,
there exists a \seq $\set{\mu_n}\sbs \Wqdual(\bdw)$ \sth $\mu_n\to
\mu$ in total variation norm. If $\mu$ is positive,
 the \seq can be chosen to be increasing.  Note that these facts imply that  $\mu$ is a trace \ifif $\abs{\mu}$
 is a trace.\1
 (iii)  $u$ is strictly moderate \ifif \ $\abs{\tr u}\in \Wqdual(\bdw)$. In this case
the boundary trace is also attained in the sense of weak convergence
in $\Wqdual(\bdw)$ of $\set{u_*^\gb:\gb\in (0,\gb_0)}$ as $\gb\to
0$. In particular \req{limtr} holds for every $\gf\in
\Wq(\bdw)\cup C(\bdw)$.\1
(iv) If $\mu:=\tr u$ and $\set{\mu_n}$
is as in (ii) then $u=\lim u_{\mu_n}$. In particular, if $u>0$
then $u$ is the limit of an  increasing \seq of strictly moderate
solutions.  \1
(v) The measure $\mu=\tr u$  is regular relative to
the $q$-topology.\1
(vi) If $u$ is positive (not necessarily strictly moderate), \req{limtr} is
valid for every $\gf\in  (\Wq_+\cap L^\infty)(\bdw)$.
\es
\Remark Assertions (i)-(iv) are well known. For proofs see \cite{MV3} which
also contains further relevant citations.

\note{Proof of (v)} If
$\mu$ is a trace then $\mu_+$ and $\mu_-$
  are traces of solutions of \rqq. Therefore it is enough to prove (v) in the case that
$\mu$ is a positive measure.

\par Every bounded Borel measure on $\bdw$ is regular in the usual sense:
$$\BAL\mu(E)=&\inf\set{\mu(O): E\sbs O,\; O \text{ relatively open}}\\
=&\sup\set{\mu(K): K\sbs E,\; K \text{ compact}} \EAL$$
for every Borel set $E\sbs \bdw$. Since
$$\BAL
\mu(E)\leq&\inf\set{\mu(D): E\sbs D\sbs \bdw,\; D \text{ \qop}}\\
\leq&\inf\set{\mu(O): E\sbs O\sbs \bdw,\; O \text{ relatively open}}
\EAL
$$
it follows that such a measure is also regular \wrto the $q$-topology.\1
\note{Proof of (vi)} By (ii) there exists an
increasing \seq of strictly moderate solutions $\set{v_n}$ \sth
$v_n\uparrow u$. If $\mu_n:=\tr v_n$ then
$$ \lim_{\gb\to 0} \int_{\Gs_\gb}v_n \gf^* dS=
  \int_{\bdw}\gf\,d\mu_n,$$
for every $\gf\in  \Wq(\bdw)$. Since $\set{\mu_n}$ increases and
converges weakly to $\mu=\tr u$ it follows that
  $\mu_n\to \mu$ in total variation. Hence
$$\int_{\bdw}\gf\,d\mu_n\to \int_{\bdw}\gf\,d\mu$$
for every bounded $\gf\in  \Wq(\bdw)$. If, in addition, $\gf\geq
0$, we obtain
$$\liminf_{\gb\to 0} \int_{\Gs_\gb}u \gf^* dS\geq
  \int_{\bdw}\gf\,d\mu.$$
On the other hand, since $u^*_\gb\to\mu$ in the sense of weak
convergence of measures, it follows that
\begin{equation}\label{tempb}
  \limsup \int_{E}u^*_\gb dS\leq \mu, \q
   \liminf \int_{A}u^*_\gb dS\geq \mu
\end{equation}
for any closed set $E\sbs \bdw$, respectively, open set
$A\sbs\bdw$). It is easily seen that, in our case, this extends to
any \qcl set $E$ (resp. \qop set $A$). Therefore if $A$ is \qop
and
$$\mu(A)=\mu(\tl A),$$
then
\begin{equation}\label{tempc}
\lim \int_{A}u_*^\gb dS = \mu(A).
\end{equation}

\par If $\gf\in  \Wq(\bdw)\cap L^\infty(\bdw)$ and $I\sbs \BBR$ is a bounded open interval then,
by \cite[Prop. 6.1.2, Prop. 6.4.10]{AH}, $A:=\gf^{-1}(I)$ is quasi
open. \Wlg we may assume that $\gf\leq 1$. Given $k\in \BBN$ and
$m=0,\ldots,2^k-1$ choose a number $a_{m,k}$ in the interval
$(m2^{-k},(m+1)2^{-k})$ \sth
$\mu\ind{\CR}(\gf^{-1}(\{a_{m,k}\})=0$. Put
$$A_{m,k}=\gf^{-1}((a_{m,k},a_{m+1,k}]) \q m=1,\ldots,2^k-1, \q A_{0,k}=\gf^{-1}((a_{0,k},a_{1,k}])$$
and
$$f_k=\sum_{m=0}^{2^k-1}m2^{-k}\chr{A_{m,k}}.$$
Then $f_k\to \gf$ uniformly and, by \req{tempc},
$$\lim_{\gb\to 0}\int_{\bdw}f_k u_*^\gb dS=\int_{\bdw}f_k d \mu\ind{\CR}.$$
This implies assertion (vi).\hfill \qedsymbol
\bth{reg-tr} Let $u\in \CU(\Gw)$.\1
\num{i} The regular set
$\CR(u)$ is \gsreg and \consy it has a regular decomposition $\set{Q_n}$. \1
\num{ii} Let
\begin{equation}\label{vCR}
  v\ind{\CR}:=\sup\set{[u]_Q:\, Q \;\text{ \qop and regular }}.
\end{equation}
Then  there exists an increasing \seq of moderate solutions $\set{w_n}$ \sth
\begin{equation}\label{vn-cond}
 \suppq w_n\sbsq \CR(u),\q w_n\uar v\ind{\CR}.
\end{equation}
 (Thus $v\ind{\CR}$ is \gsmod.)\1
\num{iii} Let $F$ be a \qcl set \sth $F\sbsq \CR(u)$. Then, for every $\ge>0$,
 there exists
a \qop regular set $Q_\ge$ \sth $\Capq(F\sms Q_\ge)<\ge$. If, in addition, $[u]_F$ is moderate then
$F$ is regular; \consy there
exists a \qop regular set $Q$ \sth $F\sbsq Q$.\1
\num{iv} With $\set{Q_n}$ as in (i), denote
\begin{equation}\label{muCRdef}
v_n:=[u]_{Q_n}, \q \mu_n:=\tr v_n, \q v:=\lim v_n,\q \muCR:=\lim \mu_n.
\end{equation}
Then,
\begin{equation}\label{v=vCR}
 v=v\ind{\CR}.
\end{equation}
Furthermore, for every \qop regular set $Q$,
\begin{equation}\label{muCRloc1}
 \mu\ind{\CR}\chr{\tl Q}=\tr [u]_Q=\tr [\vCR]_Q.
\end{equation}
Finally, $\mu\ind{\CR}$ is $q$-locally finite on $\CR(u)$ and $\gs$-finite on
$\CR_0(u):=\cup Q_n$.\1
\num{v} If $\set{w_n}$ is a \seq of moderate solutions satisfying  conditions \req{vn-cond} then,
\begin{equation}\label{muCRind}
\muCR=\lim \tr w_n
\end{equation}
\num{vi} The {\em regularized measure} $\bar\mu\ind{\CR}$ given by
\begin{equation}\label{barmuCR}
  \bar\mu\ind{\CR}(E):=\inf\set{\mu\ind{\CR}(Q): E\sbs Q,\q Q\text{ \qop}\forevery E\sbs \bdw,\; E \text{ Borel} }
\end{equation}
is $q$-perfect.\1
\num{vii}\hspace{2mm} $u\app{\CR(u)}v\ind{\CR}$.\1
\num{viii} \hspace{2mm} For every \qcl set $F\sbsq \CR(u)$:
\begin{equation}\label{consistF}
[u]_F=[v\ind{\CR}]_F.
\end{equation}
If, in addition, $\muCR(F)<\infty$ then $[u]_F$ is moderate and
\begin{equation}\label{consist-tr}
\tr[u]_F=\mu\ind{\CR}\chr{F}.
\end{equation}
\num{ix} If $F$ is a \qcl set then
\begin{equation}\label{Freg}
\muCR(F)<\infty\iff [u]_F\text{ is moderate }\iff F\text{ is regular}.
\end{equation}
\es
\bproof  (i)\hspace{2mm} By \cite[Sec. 6.5.11]{AH} the $(\ga,p)$-fine topology
possesses the quasi-Lindel\"{o}f property.  This  implies that
$\CR(u)$  is \gsreg. By \rlemma{decomposition} $\CR(u)$ has a
regular decomposition $\set{Q_n}$. Recall that $\tl Q_n\sbs Q_{n+1}$ and $\Capq(\CR(u)\sms \CR_0(u))=0$.\1
(ii)\hspace{2mm} This assertion is an immediate consequence of \req{vCR} and \rprop{sup u_n}.\1
(iii) By definition, every point in $\CR(u)$ possesses a \qop regular \ngh. Therefore, the existence of
a set $Q_\ge$,  as in the first
part of this assertion, is an immediate  consequence of \rlemma{q-compact}. Let $O_\ge$ be an open set
containing $F\sms Q_\ge$ \sth $\Capq(O_\ge)<2\ge$. Put $F_\ge:=F\sms O_\ge$. Then $F_\ge$ is a
\qcl set,
$F_\ge\sbs F$, $\Capq(F\sms F_\ge)<2\ge$ and $F_\ge\sbsq Q_\ge$.

\note{Assertion 1} {\em Let $E$ be a \qcl set, $D$ a \qop regular set and $E\sbsq D$. Then there exists a
decreasing \seq of \qop sets $\set{G_{n}}_{n=1}^\infty$ \sth
\begin{equation}\label{regtrA}
 E\sbsq G_{n+1}\sbs \tl G_{n+1}\sbsq G_{ n}\sbsq D
\end{equation}
and
\begin{equation}\label{regtrB}
  [u]_{G_{n}}\to [u]_{E} \txt{in} L^q(\Gw,\gr).
\end{equation}
}
\par By \rlemma{set-approx} and \rth{lim[u]}, there exists a decreasing \seq of \qop sets
$\set{G_{n}}$ satisfying \req{regtrA} and, in addition,
\sth $[u]_{G_{n}}\dar [u]_{E}$ locally uniformly in $\Gw$. Since
$[u]_{G_{n}}\leq [u]_{D}$ and the latter is a moderate solution
we obtain \req{regtrB}.

Put
$$E_n:=\cup_{m=1}^n F_{1/m}, \q D_n:=\cup_{m=1}^n Q_{1/m}.$$
Then $E_n$ is \qcl, $D_n$ is \qop and regular and $E_n\sbsq D_n$. Therefore, by Assertion 1,
it is possible to
choose a \seq of \qop regular sets $\set{V_n}$ \sth
\begin{equation}\label{regtrC}
E_n\sbsq V_n\sbs \tl V_n\sbsq D_n,\q
\norm{[u]_{V_n}-[u]_{E_n}}\ind{L^q(\Gw,\gr)}\leq 2^{-n}.
\end{equation}
By \rth{lim[u]},
$$ [u]_F\leq [u]_{E_n}+[u]_{F\sms E_n} \txt{and} [u]_{F\sms E_n}\dar 0.$$
Therefore
$[u]_{E_n}\uar [u]_F$.

If, in addition, $[u]_F$ is moderate then
$$[u]_{E_n}\uar [u]_F \txt{in} L^q(\Gw,\gr)$$
and \consy, by \req{regtrC}, $$[u]_{V_n}\to [u]_F \txt{in} L^q(\Gw,\gr).$$
 Let $\set{V_{n_k}}$
be a \sseq \sth
\begin{equation}\label{regtrD}
  \norm{[u]_{V_{n_k}}- [u]_F}_{L^q(\Gw,\gr)}\leq 2^{-k}.
\end{equation}
Recall that $E_n\sbsq V_n\cap F$ and that $\Capq(F\sms E_n)\dar 0$.
Therefore $\Capq(F\sms V_{n})\to 0$. \Consy
$F\sbsq W:=\cap_{k=1}^\infty V_{n_k}$ and, in view of \req{regtrD}, $[u]_W $ is moderate.
Obviously this implies that $W$ is pre-regular (any
\qcl set $E\sbsq W$ has the property that $[u]_E$ is moderate) and $F$ is regular. Finally, by
\rlemma{e-ngh}, every \qcl regular set is contained in a \qop regular set.\1
(iv)\hspace{2mm}  Let $Q$ be a \qop regular set and put
$\mu\ind{Q}=\tr [u]_Q$. If $F$ is a \qcl set \sth $F\sbsq Q$ then, by \rth{maintrace},
\begin{equation}\label{regtr1}
 [u]_{ F}=\mu\ind{Q}\chr{F}.
\end{equation}
In particular the compatibility condition holds: if $Q,Q'$ are \qop regular sets then
\begin{equation}\label{regtr3}
 \mu_{Q\cap Q'}=\mu_{Q}\chr{\tl Q\cap \tl Q'}=\mu_{Q'}\chr{\tl Q\cap \tl Q'}.
\end{equation}
With the notation of \req{muCRdef}, $[v_{n+k}]_{Q_k}=v_k$ and hence  $\mu_{n+k}\chr{\tl Q_k}=\mu_k$ for every $k\in\BBN$.
\par Let $F$ be an arbitrary \qcl subset of $\CR(u)$.  Since $\Capq (F\sms Q_n)\to 0$ it follows that
\begin{equation}\label{regtrE}
 [v_n]_F=[u]_{F\cap \tl Q_n}\uar [u]_F.
\end{equation}
In addition, $[v]_F\geq \lim [v_n]_F=[u]_F$ and $v\leq u$ lead to,
\begin{equation}\label{regtrF}
 [u]_F=[v]_F.
\end{equation}
If $Q$ is a \qop regular set,
$[u]_{Q}=\lim[v_n]_Q\leq \lim v_n=:v$ and so $v\ind{\CR}\leq v$. On the other hand it is
obvious that $v\leq v\ind{\CR}$. Thus \req{v=vCR} holds.

By \req{regtr1} and \req{regtrE}, if $F$ is a \qcl subset of $\CR(u)$ and $[u]_F$ is moderate,
\begin{equation}\label{regtrG}
  \tr [u]_F=\lim\tr[v_n]_F=\lim\mu_n\chr{F}= \mu\ind{\CR}\chr{F},
\end{equation}
which implies \req{muCRloc1}.  This also shows that
$\mu\ind{\CR}\chr{F}$ is independent of the choice of the \seq $\set{\mu_n}$ used in its definition.
This remains valid
for any \qcl set  $ F\sbsq \CR(u)$ because $\Capq(\CR(u)\sms \CR_0(u))=0$ and $\mu\ind{\CR}$ is $\gs$-finite on $\CR_0(u)$.
The last observation is a consequence of  the fact that $\CR_0(u)$
has a regular decomposition.

\par The assertion that $\mu\ind{\CR}$ is
$q$-locally finite on $\CR(u)$ is a consequence of the fact that
every point in $\CR(u)$ is contained in a \qop regular set.  \1
(v) If $w$ is a moderate solution and $w\leq \vCR$ and $\suppq w\sbsq \CR(u)$
then $\tau:=\tr w\leq \muCR$. Indeed
$$[w]_{Q_n}\leq [\vCR]_{Q_n}=v_n,\; [w]_{Q_n}\uar w \Lra \tr[w]_{Q_n}\uar \tau\leq \lim\tr v_n=\muCR.$$

Now, let $\set{w_n}$ be an increasing \seq of moderate solutions \sth $F_n:=\suppq w_n\sbsq \CR(u)$
and $w_n\uar v\ind{\CR}$. We must show that, if $\nu_n:=\tr w_n$,
\begin{equation}\label{nu=mu}
  \nu:=\lim \nu_n=\muCR.
\end{equation}
\par By the previous argument  $\nu\leq \mu\ind{\CR}$.
The opposite inequality is obtained as follows.
Let $D$ be a \qop regular set and let $K$ be a compact subset of $D$ \sth $\Capq(K)>0$.
$$w_n\leq [w_n]_D+[w_n]_{D^c} \Lra \vCR=\lim w_n\leq \lim[w_n]_D+ U_{D^c}.$$
The \seq $\set{[w_n]_D}$ is dominated by the moderate function $[\vCR]_D$. In addition
$\tr[w_n]_D=\nu_n\chr{\tl D}\uar \nu\chr{\tl D}$. Hence, $\nu\chr{\tl D}$ is a bounded measure and
 $[w_n]_D\uar u_{\nu\chr{\tl D}}$ where the function on the right is the moderate solution with trace
$\nu\chr{\tl D}$. \Consy
$$\vCR=\lim w_n\leq u_{\nu\chr{\tl D}}+ U_{D^c}.$$
This in turn implies
$$([\vCR]_K- u_{\nu\chr{\tl D}})_+ \leq\inf(U_{D^c},U_K)$$
the function on the left being a subsolution and the one on the left a supersolution. Therefore
$$([\vCR]_K- u_{\nu\chr{\tl D}})_+\leq [[U]_{D^c}]_K=0.$$
Thus, $[\vCR]_K\leq u_{\nu\chr{\tl D}}$ and hence $\muCR\chr{K}\leq \nu\chr{\tl D}$.
Further, if $Q$ is a \qop set \sth $\tl Q\sbsq D$ then, in view of the fact that
$$\sup\set{\muCR\chr{K}: K\in Q,\;K\text{ compact}}=\mu\chr{Q},$$
we obtain,
\begin{equation}\label{regtrH}
\muCR\chr{Q}\leq \nu\chr{\tl D}.
\end{equation}
Applying this inequality to the sets $Q_m, Q_{m+1}$ we finally obtain
$$\muCR\chr{Q_m}\leq \nu\chr{\tl Q_{m+1}}\leq \nu\chr{Q_{m+2}}.$$
Letting $m\tin$ we conclude that $\muCR\leq \nu$. This completes the proof of \req{nu=mu} and of
assertion (v).\1
(vi) The measure $\mu\ind{\CR}$ is essentially absolutely continuous relative to $\Capq$ (see \rdef{cap-reg}).
Therefore this assertion follows from \rlemma{cap-reg}.\1
(vii)  By \req{split2}
$$ u\leq [u]_{Q_n}+[u]_{\bdw\sms  Q_n}.$$
 By \rth{sup/inf}(c)
 $$ [u]_{\bdw\sms Q_n}\dar [u]_{\bdw\sms \CR_0}.$$
 Hence
 $$ \lim(u-[u]_{Q_n})= u-v\ind{\CR}\leq [u]_{\bdw\sms \CR_0}$$
 so that $u\ominus v\ind{\CR}\app{\CR_0}0$. Since $\vCR\leq u$ this is equivalent to the statement
$u\app{\CR}\vCR$.\1
(viii) \req{consistF} was established before, see \req{regtrF}. Alternatively, it follows from the previous assertion and \rth{uappv}.
If $\mu\ind{\CR}(F)$ is finite then  \req{consist-tr} is a consequence of (i) and \req{regtr1}. Indeed, $F_n:=F\cap \tl Q_n\uar F_0\smq F$.
Hence, $[u]_F\leq [u]_{F_n}+[u]_{F\sms F_n}$ and $\Capq(F\sms F_n)\dar 0$. Hence $[u]_F=\lim[u]_{F_n}$ and
$\tr[u]_{F_n}=\muCR\chr{F_n}\uar \muCR\chr{F_0}=\muCR\chr{F}$. Since $\muCR\chr{F}$
is a bounded measure, $[u]_F$ is moderate and \req{consist-tr} holds.\1
(ix) If $\muCR(F)<\infty$ then, by (viii), $[u]_F$ is moderate and, by (iii), $F$ is regular.
Conversely, if $F$ is regular then $[u]_F$ is moderate and, by \req{muCRloc1}, $\muCR(F)<\infty$.
\end{proof}
\subsection{The precise boundary trace}
\bdef{precisetr}{\rm
 Let $q_c\leq q$ and $u\in \CU(\Gw)$.
  \begin{description}
\item[a]  The solution $v\ind{\CR}$ defined by \req{vCR} is called the {\em regular component} of
$u$ and will be denoted by $u\ind{reg}$.
\item[b] Let $\set{v_n}$ be an increasing \seq of moderate solutions satisfying condition
\req{vn-cond} and put $\mu\ind{\CR}=\muCR(u):=\lim \tr v_n$. Then,
the regularized measure $\bar\mu\ind{\CR}$,
  defined by \req{barmuCR}, is called the
  {\em regular boundary trace} of $u$. It will be denoted by $\trR u$.
\item[c] The couple $(\trR u,\CS(u))$
is called the {\em precise boundary trace} of $u$ and will be denoted by
$\tr^c u$.
\item[d] Let $\nu$ be the Borel measure on $\bdw$ given by
\begin{equation}\label{measure-tr}
\nu(E)=\begin{cases}(\tr\ind{\CR}u)(E)& \text{if }\;E\sbs\CR(u),\\
\nu(E)=\infty &\text{if }\;E\cap \CS(u)\neq\ems,
\end{cases}
\end{equation}
for every Borel set $E\sbs \bdw$. Then $\nu$ is the measure representation of the
precise boundary trace of $u$, to be denoted by $\tr u$.
\end{description}
Note that, by \rth{reg-tr} (v),
  the measure $\mu\ind{\CR}$ is independent of the choice of the \seq $\set{v_n}$.
}\es
\bth{gsmod-tr} Assume that $u\in \CU(\Gw)$ is a \gsmod solution, i.e.,
there exists an increasing \seq $\set{u_n}$ of positive moderate
solutions \sth $u_n\to u$. Let $\mu_n:=\tr u_n$, $\mu_0:=\lim
\mu_n$ and put
\begin{equation}\label{modreg}
\mu(E):=\inf\set{\mu_0(D): E\sbs D\sbs \bdw,\; D \text{ \qop}},
\end{equation}
for every Borel set $E\sbs \bdw$. Then: \1
\num{i} $\mu$ is the precise boundary trace of $u$ and $\mu$ is $q$-perfect.
In particular $\mu$ is independent of the \seq $\set{u_n}$ which appears in its definition.\1
\num{ii} If $A$ is a Borel set \sth $\mu(A)<\infty$ then $\mu(A)=\mu_0(A)$.\1
\num{iii} A solution $u\in \CU(\Gw)$ is  \gsmod \ifif
\begin{equation}\label{modreg1}
 u=\sup\set{v\in \CU(\Gw): \text{$v$ moderate}\q v\leq u},
\end{equation}
which is equivalent to
\begin{equation}\label{modreg2}
u=\sup\set{u_\tau:\tau\in \Wqdb,\q \tau\leq \tr u}.
\end{equation}
Thus, if $u$ is \gsmod, there exists an increasing \seq of {\em strictly
moderate} solutions converging to $u$.\1
\num{iv} If $u,w$ are \gsmod solutions,
\begin{equation}\label{smod-monotone}
 \tr w\leq \tr u \iff w\leq u.
\end{equation}
\es
\bproof  (i)\hspace{2mm} Let $Q$ be a \qop set and $A$ a Borel set
\sth $\Capq(A)=0$.
Then $\mu_n( A)=0$ so that $\mu_0(A)=0$. Therefore $\mu_0(Q\sms A)=\mu_0(Q)$. Thus $\mu_0$ is essentially
\abc and, by  \rlemma{cap-reg}, $\mu$ is $q$-perfect.
\par Let  $\set{D_n}$ be a regular decomposition of
$\CR(u)$. Put $D'_n=\CR(u)\sms D_n$ and observe that $D'_n\dar E$ where $\Capq(E)=0$. Therefore
$$u_{\mu_n\chr{D'_n}}\dar 0.$$
Since,
$$\mu_n\chr{\CR(u)}=\mu_n\chr{D_n}+ \mu_n\chr{D'_n},$$
it follows that
$$ \abs{u_{\mu_n\chr{\CR(u)}}- \lim u_{\mu_n\chr{D_n}}}\leq u_{\mu_n\chr{D'_n}}\to 0.$$
Now
$$u_n\leq u_{\mu_n\chr{\CR(u)}}+[u]_{\CS(u)}.$$
Hence
$$u-[u]_{\CS(u)}\leq w:=\lim u_{\mu_n\chr{\CR(u)}}=\lim u_{\mu_n\chr{D_n}}\leq u\ind{reg}.$$
This implies $u\ominus[u]_{\CS(u)}\leq u\ind{reg}$. But the definition of $u\ind{reg}(=v\ind{\CR})$ implies that
  $u\ind{reg}\leq u\ominus[u]_{\CS(u)}$. Therefore
  $\lim u_{\mu_n\chr{D_n}}=u\ind{reg}$.
Thus the \seq $\set{u_{\mu_n\chr{D_n}}}$ satisfies condition \req{vn-cond} and \consy, by \rth{reg-tr} (iv) and \rdef{precisetr},
\begin{equation}\label{gsmod1}
 \lim \mu_n\chr{D_n}=\mu\ind{\CR},\q \tr\!\ind{\CR}u=\bar\mu\ind{\CR}.
\end{equation}

\par If $\gx\in\CS(u)$ then, for every \qop \ngh $Q$ of $\gx$,
$\int_{\Gw}[u]_Q^q\gr dx=\infty$. This implies: $\mu_n(\tl Q)\tin$. To verify this fact, assume that, on the contrary,
 there exists
a \sseq (still denoted $\set{\mu_n}$) \sth $\sup \mu_n(\tl Q)<\infty$. Denote by $v_{n,Q}$ and $w_{n,Q}$
 the solutions with boundary trace $\chr{Q}\tr u_n$
and $\chr{Q^c}\tr u_n$ respectively. Then
$$u_n\leq v_{n,Q}+w_{n,Q}\Lra u\leq v_Q+w_Q,\q v_Q=\lim v_{n,Q}, \q w_Q=\lim w_{n,Q}.$$
 Then $v_Q$ is moderate and  $w_Q$ vanishes
in $Q$. If $D$ is a \qop set \sth $\gx\in D\sbs \tl D\sbs Q$ then $[w_Q]_D=0$. Therefore
$$\min(u,U_D)\leq \min (v_{Q},U_{D})+ \min (w_{Q},U_{D})\Lra [u]_D\leq [v_Q]_D,$$
which brings us to a contradiction.
In conclusion,  if $\gx\in\CS(u)$ then $\mu_0(\tl Q)=\infty$ for every \qop \ngh of $\gx$. \Consy
$\mu({\gx})=\infty$. This fact and \req{gsmod1} imply that $\mu$ is the precise trace of $u$.\1
(ii)\hspace{2mm} If $\mu(A)<\infty$ then $A$ is contained in a \qop set $D$ \sth $\mu_0(D)<\infty$ and,
by \rlemma{cap-reg}, $\mu(A)=\mu_0(A)$.\1
(iii)\hspace{2mm} Let $u\in \CU(\Gw)$ be \gsmod and put
\begin{equation}\label{u*}
 u^*:=\sup\set{v\in \CU(\Gw): \text{$v$ moderate}\q v\leq u}.
\end{equation}
By its definition $u^*\leq u$. On the other hand, since there exists an increasing \seq of
moderate solutions $\set{u_n}$ converging to $u$, it follows
that $u\leq u^*$. Thus $u=u^*$.

Conversely, if $u\in \CU(\Gw)$ and $u=u^*$ then, by \rprop{sup u_n},
there exists an increasing \seq of moderate solutions $\set{u_n}$ converging to $u$. Therefore $u$ is \gsmod.

In view of \rprop{mod-tr} (iv), 
$$u^*\leq \sup\set{u_\tau:\tau\in \Wqdb,\q \tau\leq \tr u}=:u^\ddag.$$
On the other hand, if $u$ is \gsmod, $\tau\in \Wqdb$ and $\tau\leq \tr u$ then
(with $\mu_n$ and $u_n$ as in the statement of the theorem), $\tr (u_\tau\ominus u_n)=(\tau-\mu_n)_+\dar 0$.
Hence $u_\tau\ominus u_n\dar 0$, which implies that $u_\tau\ominus u=0$, i.e. $u_\tau\leq u$. Therefore
$u\geq u^\ddag.$
Thus \req{modreg1} implies \req{modreg2} and each of them implies that $u$ is \gsmod. Therefore the two are equivalent.\1
(iv)\hspace{2mm} The assertion $\Lra$ is a consequence of \req{modreg2}. To verify the assertion $\Lla$
 it is sufficient to show that if  $w$ is moderate, $u$ is \gsmod and $w\leq u$ then $\tr w\leq \tr u$.
Let $\set{u_n}$ be an increasing \seq of positive moderate solutions converging to $u$. Then
$u_n\vee w\leq u$  and \consy  $u_n\leq u_n\vee w\uar u$.
Therefore $\tr(u_n\vee w)\uar \mu'\leq \tr u$ so that $\tr w\leq \tr u$.

\eproof

\bth{precisetr} Let $u\in \CU(\Gw)$ and put $\nu=\tr u$.\1
\num{i} $u\ind{reg}$ is \gsmod and $\tr u\ind{reg}=\trR u$.\1
 \num{ii} If $v\in \CU(\Gw)$
\begin{equation}\label{precise1}
 v\leq u \Lra
\tr v\leq \tr u.
\end{equation}
 If $F$ is a \qcl set then
\begin{equation}\label{precise2}
\tr[u]_F\leq \nu\chr{F}.
\end{equation}
\num{iii} A singular point can be characterized in terms of the measure $\nu$ as follows:
\begin{equation}\label{precise3}
 \gx\in\CS(u)\iff \nu(Q)=\infty \forevery Q: \; \gx\in Q,\; Q\text{ \qop}.
\end{equation}
\num{iv}  If $Q$ is a \qop set then:
\begin{align}\label{precise4}
Q\;\text{\rm pre-regular}&\iff \nu(F\cap Q)<\infty \forevery F\sbsq Q: \;F \;\text{\rm\qcl}\\
Q\;\text{\rm regular}&\iff \exists\; \text{\rm Borel set $A$: }\Capq(A)=0,\;\nu(\tl Q\sms A)<\infty.\label{precise5}
\end{align}
\num{v} The singular set of $u\ind{reg}$ may not be empty. In fact
\begin{equation}\label{precise6}
 \CS(u)\sms b_q(\CS(u))\sbs\CS(u\ind{reg})\sbs \CS(u)\cap\wtl{\CR(u)},
\end{equation}
where $b_q(\CS(u))$ is the set of $\Capq$-thick points of $\CS(u)$, (see Notation 2.1). \1
\num{vi} Put
\begin{equation}\label{precise7}
 \CS_0(u):=\set{\gx\in\bdw:\, \nu(Q\sms\CS(u))=\infty \forevery Q: \; \gx\in Q,\; Q\text{ \qop}}.
\end{equation}
Then
\begin{equation}\label{precise8}
\CS(u\ind{reg})\sms b_q(\CS(u))\sbs\CS_0(u)\sbs \CS(u\ind{reg})\cup b_q(\CS(u))
\end{equation}
\es
\Remark This result complements \rth{reg-tr} which deals with the regular boundary trace.
\bproof
\num{i} By \rth{reg-tr} (ii) $u\ind{reg}$ is \gsmod. The second part of the statement follows from
\rdef{precisetr} and \rth{gsmod-tr} (i).\1
\num{ii} If $v\leq u$ then $\CR(u)\sbs \CR(v)$ and by definition $v\ind{reg}\leq u\ind{reg}$. By \rth{gsmod-tr} (iv)
$\tr v\ind{reg}\leq \tr u\ind{reg}$  and \consy $\tr v\leq \tr u$.
 \req{precise2} is an immediate consequence of \req{precise1}.\1
 \num{iii} If $\gx\in \CR(u)$ there exists a \qop regular \ngh $Q$ of $\gx$. Hence
 $\nu(Q)=\trR u(Q)<\infty$. If $\gx\in \CS(u)$, it follows immediately from the
 definition of precise trace that $\nu(Q)=\infty$ for every \qop \ngh $Q$ of $\gx$.\1
\num{iv} If $Q$ is pre regular then $[u]_F$ is moderate for
every \qcl set $F\sbsq Q$ and $Q\sbsq \CR(u)$. By \rth{reg-tr} (viii) this implies:
$\tr [u]_F=(\muCR(u))\chr{F}$ and \consy $\nu(F\cap Q)=(\muCR(u))(F\cap Q)=(\muCR(u))(F)<\infty$.
 Therefore \req{precise4} holds in the direction $\Lra$.

  Conversely, if $Q$ is a \qop set, $F$ a \qcl set, $F\sbsq Q$ and $\nu(F\cap Q)<\infty$ then,
  by \rdef{precisetr}, $F\cap Q\sbs \CR(u)$
  which implies $F\sbsq \CR(u)$ and $\muCR(u)(F)<\infty$. Therefore, by \rth{reg-tr} (viii)
$[u]_F$ is moderate. This implies \req{precise4} in the opposite direction.

If $Q$ is regular there exists a pre-regular set $D$ \sth $\tl Q\sbsq D$. Therefore \req{precise4} implies
\req{precise5} in the direction $\Lra$. On the other hand,
$$\nu(\tl Q\sms A)<\infty\Lra \tl Q\sbsq \CR(u)$$
and $\muCR(\tl Q)=\muCR(\tl Q\sms A)<\infty$. Hence, by \rth{reg-tr} (ix), $[u]_Q$ is moderate and
 $\tl Q$ is regular.\1
\num{v} Since $\suppq u\ind{reg}\sbs \wtl{\CR(u)}$ and $\CR(u)\sbs\CR(u\ind{reg})$ we have
$$\CS(u\ind{reg})\sbs \wtl{\CR(u)}\cap \CS(u).$$
Next we show that $\CS(u)\sms b_q(\CS(u))\sbs\CS(u\ind{reg})$.

If $\gx\in \CS(u)\sms b_q(\CS(u))$
then  $\CR(u)\cup\set{\gx}$ is a \qop \ngh of $\gx$.
By (i) $u\ind{reg}$ is \gsmod and \consy
(by \rth{gsmod-tr} (i)) its trace is $q$-perfect.
Therefore, if $Q_0$ is a \qop \ngh of $\gx$ and $Q=Q_0\cap(\set{\gx}\cup\CR(u))$ then
$$(\tr u\ind{reg})(Q)=(\tr u\ind{reg})(Q\sms \set{\gx})=(\tr u)(Q\sms \set{\gx}).$$
The last equality is valid because $Q\sms \set{\gx}\sbs \CR(u)$. Let $D$ be a \qop set
\sth $\gx\in D\sbs \tl D\sbs Q$. If $\tr u(\tl D\sms \set{\gx})<\infty$ then, by (iv),
$D$ is regular and $\gx\in \CR(u)$, contrary to our assumption. Therefore $\tr u(\tl D\sms \set{\gx})=\infty$
so that $\tr u\ind{reg}(Q_0\sms \set{\gx})=\infty$ for every \qop \ngh $Q_0$ of $\gx$,
 which implies  $\gx\in \CS(u\ind{reg})$. This completes the proof of \req{precise6}.\1
\num{vi}
 If $\gx\nin b_q(\CS(u))$, there exists a \qop
  \ngh $D$ of $\gx$ \sth $(D\sms\set{\gx})\cap\CS(u)=\ems$ and \consy
\begin{equation}\label{precise9}
 (\tr u\ind{reg})( D\sms\set{\gx})=(\tr u\ind{reg})( D\sms\CS(u))=(\tr u)( D\sms\CS(u)).
\end{equation}
If, in addition, $\gx\in \CS_0(u)$ then
$$(\tr u)( D\sms\CS(u))=(\tr u\ind{reg})( D\sms\set{\gx})=\infty.$$
  If $Q$ is an arbitrary  \qop \ngh  of $\gx$ then the same holds if $D$ is replaced by $Q\cap D$.
  Therefore $(\tr u\ind{reg})( Q\sms\set{\gx})=\infty$ for any such $Q$.
  \Consy $\gx\in\CS(u\ind{reg})$, which proves that $\CS_0(u)\sms b_q(\CS(u))\sbs \CS(u\ind{reg})$.

  On the other hand, if $\gx\in \CS(u\ind{reg})\sms b_q(\CS(u))$ then there exists a \qop \ngh
  $D$ \sth \req{precise9} holds and $(\tr u\ind{reg})(D)=\infty$.
Since $u\ind{reg}$ is \gsmod, $(\tr u\ind{reg})$ is $q$-perfect so that
$(\tr u\ind{reg})(D)=(\tr u\ind{reg})(D\sms\set{\gx})=\infty$. \Consy, by \req{precise9},
$(\tr u)( D\sms\CS(u))=\infty$. If $Q$ is any \qop \ngh of $\gx$ then $D$ can be replaced by
$D\cap Q$. \Consy $(\tr u)(Q\sms\CS(u))=\infty$ and we conclude that $\gx\in \CS_0(u)$.
This completes the proof of \req{precise8}.

\eproof


\blemma{bqF}
Let $F\sbs\bdw$ be a \qcl set. Then $\CS(U_F)=b_q(F)$.
\es
\bproof Let $\gx$ be a point on $\bdw$ \sth $F$ is $\Capq$-thin at $\gx$. Let $Q$ be a \qop \ngh
of $\gx$ \sth $\tl Q\sbsq F^c$. Then $[U_F]_{Q}=U_{F\cap \tl Q}=0$. Therefore $\gx\in \CR(U_F)$.
\par Conversely, assume that $\gx\in F\cap\CR(U_F)$. By \rth{loc-tr2}
 there exists a
\qop \ngh $Q$ of $\gx$ \sth
 $u^Q=\lim_{\gb\to 0} u_\gb^Q$ exists and $u^Q$ is a moderate solution. Let $D$ be a \qop \ngh
 of $\gx$ \sth $\tl D\sbsq Q$. Then \rlemma{convergence2} implies that $[u]_D$ is moderate so that
 $D\sbs \CR(u)$. In turn this implies that $\Capq(F\cap D)=0$ and \consy $F$ is $q$-thin at $\gx$.
\eproof
\subsection{The boundary value problem}\ \1
\note{Notation 5.2}
\begin{description}
\item[{\rm a}] Denote by $\GTM_+(\bdw)$  the space of positive Borel
measures on $\bdw$ (not necessarily bounded).
\item[{\rm b}] Denote by $\GTC_q(\bdw)$ the space of couples $(\tau,F)$ \sth
 $F\sbs\bdw$ is  \qcl,
$\tau\in \GTM_+(\bdw)$, $\qsupp \tau\sbs\wtl{\bdw\sms F}$ and $\tau\chr{\bdw\sms F}$
is  $q$-locally finite.
\item[{\rm c}] Denote by $\BBT:\GTC_q(\bdw)\mapsto \GTM_+(\bdw)$
 the mapping given by $\nu=\BBT(\tau, F)$ where $\nu$ is defined as in \req{measure-tr} with
  $\CS(u),\CR(u)$ replaced by $F, F^c$ respectively. $\nu$ is the measure representation of the
  couple $(\tau, F)$.
\item[{\rm d}] If $(\tau, F)\in \GTC_q(\bdw)$ the set
\begin{equation}\label{Ftau}
 F_\tau=\set{\gx\in \bdw: \tau(Q\sms F)=\infty \q\forall Q\txt{\qop \ngh of $\gx$}\!}
\end{equation}
is called the set of explosion points of $\tau$.
\end{description}
\Remark Note that $F_\tau\sbs F$ (because $\tau\chr{\bdw\sms F}$
is  $q$-locally finite) and $F_\tau\sbs\wtl{\bdw\sms F}$ (because $\tau$ vanishes outside this set).
Thus
\begin{equation}\label{Ftau-relation}
 F_\tau\sbs b_q(\bdw\sms F)\cap F.
\end{equation}

\bth{existence} Let $\nu$ be a positive Borel measure on $\bdw$.\1
\num{i} The boundary value problem
\begin{equation}\label{generalBVP}
 -\Gd u + u^q=0,\q u>0 \text{ in } \Gw,\q \tr (u)=\nu \text{ on }\bdw
\end{equation}
possesses a solution \ifif $\nu\in \BBM_q(\bdw)$. \1
\num{ii} Let $(\tau,F)\in \GTC_q(\bdw)$ and put
$\nu:=\BBT(\tau,F)$. Then $\nu\in \BBM_q(\bdw)$ \ifif
\begin{equation}\label{nesc-exist}
  \tau\in \BBM_q(\bdw),\q F=b_q(F)\cup F_\tau.
\end{equation}
\num{iii} Let $\nu\in \BBM_q(\bdw)$ and denote
\begin{equation}\label{CFnu}\BAL
 \CE_\nu&:=\set{E: E \text{ $q$-quasi-closed, }\; \nu(E)<\infty},\\
 \CD_\nu &:=\set{D : \, D\text{ \qop }, \; \tl D\smq E \text{ for some }E\in \CE_q}.
\EAL\end{equation}
Then a solution of \req{generalBVP}
is given by $u=v\oplus U_{F}$ where
\begin{equation}\label{existence1}
  G:=\bigcup_{\CD_\nu}D,\q F:=\bdw\sms G,\q v:=\sup\set{u_{\nu \chr{E}}: \, E \in \CE_\nu}.
\end{equation}

Note that if $E\in \CE_\nu$ then $\nu \chr{E}$ is a bounded Borel measure which
does not charge sets of $\Capq$-capacity zero. Recall that if $\mu$ is a positive measure
possessing these properties then $u_\mu$ denotes the moderate solution with boundary trace $\mu$.\1
\num{iv} The solution $u=v\oplus U_{F}$ is \gsmod and it is the unique solution of problem
\req{generalBVP} in the class of \gsmod solutions. Furthermore, $u$ is the largest
solution of the problem.
\bcom
\begin{equation}\label{existence1'}\BAL
\CF_\nu&=\set{Q : \, Q\textrm{ \qop }, \;\nu(Q)<\infty},& G&=\bigcup_{\CF_\nu}Q,\\
 v&=\sup\set{u_\mu: \mu=\nu \chr{Q},\; Q\in \CF_\nu},& F&=\bdw\sms G.
\EAL\end{equation}
\end{comment}
\es
\bproof  First we prove\1 
\num{\bf A}{\em If $u\in \CU(\Gw)$}
\begin{equation}\label{(i)onlyif}
 \tr u=\nu \; \Lra \;\nu\in \BBM_q(\bdw).
\end{equation}

By \rth{reg-tr}, $u\ind{reg}$ is \gsmod and $u\app{\CR(u)}u\ind{reg}$. Therefore
$$(\tr u)\chr{\CR(u)}= (\tr u\ind{reg})\chr{\CR(u)}.$$
By \rth{gsmod-tr},  $\bar\mu\ind{\CR}:=\tr u\ind{reg}\in\BBM_q(\bdw)$. If $v$ is defined
as in \req{existence1} then
\begin{equation}\label{existA}
  v=\sup\set{[u]_Q: Q\text{ \qop regular set}}=u\ind{reg},
\end{equation}
where the second equality holds by definition.
Indeed, by \rth{precisetr}, for every \qop set $Q$,
 $[u]_Q$ is moderate \ifif $\nu(\tl Q\sms A)<\infty$ for some set $A$ of capacity zero.
 This means that
 $[u]_Q$ is moderate \ifif there exists $E\in\CE_\nu$ \sth $Q\smq E$. When this is the case,
 $$\tr [u]_Q=\muCR(u)\chr{\tl Q}=\muCR(u)\chr{E}=\nu\chr{E}.$$
Thus $v\geq u\ind{reg}$. On the other hand,
if $E\in \CE_\nu$, then $\tl E\sbsq \CR(u)$ and $\muCR(u)(\tl E)=\muCR(u)(E)<\infty$.
Therefore, by \rth{reg-tr} (ix), $\tl E$ is regular, i.e.,
there exists a \qop regular set $Q$ \sth $E\sbsq Q$. Hence $u_{\nu\chr{E}}\leq [u]_Q$ and we conclude that
$v\leq u\ind{reg}$. This proves \req{existA}.
In addition, if $E\cap\CS(u)\neq\ems$ then, by \rdef{precisetr},
$\nu(E)=\infty$. Therefore $\nu$ is outer regular \wrto  $q$-topology.
\par Next we must show that $\nu$ is essentially \abc. Let $Q$ be a \qop set and $A$ a non-empty Borel subset of $Q$ \sth $\Capq(A)=0$.
Either $\nu(Q\sms A)=\infty$ in which case
$\nu(Q\sms A)=\nu(Q)$ or $\nu(Q\sms A)<\infty$.
In the second case $Q\sms A\sbs \CR(u)$ and
$$\nu(Q\sms A)=\bar\mu\ind{\CR}(Q\sms A)=\bar\mu\ind{\CR}(Q)<\infty.$$
Let $\gx\in A$ and let $D$ be a \qop subset of $Q$ \sth $\gx\sbs D\sbs \tl D\sbs Q$. Let $B_n$ be a \qop \ngh of
$A\cap \tl D$ \sth $\Capq(B_n)<2^{-n}$ and $B_n\sbsq D$. Then
$$[u]_D\leq [u]_{E_n}+[u]_{B_n},\qq E_n:=\tl D\sms B_n.$$
Since $\lim [u]_{B_n}=0$ it follows that $[u]_D=\lim [u]_{E_n}$. But
$$\norm{[u]_{E_n}}_{L^q(\Gw,\gr)}\leq C\nu(E_n) \leq C\nu(Q\sms A).$$
By assumption, $\nu(Q\sms A)<\infty$. Therefore $\norm{[u]_{D}}_{L^q(\Gw,\gr)}<\infty$ which
implies
that $D\sbs \CR(u)$. Since every point of $A$ has a \ngh $D$ as above we conclude that $A\sbs \CR(u)$
and  hence $\nu(A)=(\trR u)(A)=0$.
In conclusion $\nu$ is essentially \abc and $\nu\in \BBM_q(\bdw)$.

Secondly we prove:\1
\num{\bf B}{\em Suppose that  $(\tau,F)\in \GTC_q(\bdw)$ satisfies \req{nesc-exist} and put $\nu:=\BBT(\tau,F)$.
Then the solution
$u=v\oplus U_F$, with $v$ as in \req{existence1}, satisfies $\tr u=\nu$.
\par By \req{(i)onlyif}, this also implies that $\nu\in \BBM_q(\bdw)$.}\2
\indent Clearly $v$ is a \gsmod solution. The fact that $\tau$ is $q$-locally finite in
 $F^c$ and  essentially \abc relative to $\Capq$ implies that
\begin{equation}\label{exist3}
  G:=\bdw\sms F\sbs\CR(v),\q (\tr v)\chr{G}=\tau\chr{G}.
\end{equation}
It follows from the definition of $v$ that $F_\tau\sbs \CS(v)$.
Hence, by \rlemma{bqF},
\begin{equation}\label{exist5}
 F_\tau\cup b_q(F)\sbs \CS(v)\cup \CS(U_F)\sbs \CS(u)\sbs F.
\end{equation}
Hence, by \req{nesc-exist}, $\CS(u)=F$, $v=u\ind{reg}$ and $\tau=\tr u\ind{reg}$.
Thus $\tr^c u=(\tau,F)$
which is equivalent to $\tr u=\nu$.
\par Next we show:\1
\num{\bf C}{\em Suppose that $(\tau,F)\in \GTC_q(\bdw)$ and that there exists a solution $u$ \sth
$\tr^cu=(\tau,F)$. Then,
\begin{equation}\label{exist6}
 \tau=\trR u=\tr u\ind{reg},\q F=\CS(u).
\end{equation}
If $U:=u\ind{reg}\oplus U_F$ then $\tr U=\tr u$ and $u\leq U$. $U$ is the unique \gsmod solution
of \req{generalBVP} and $(\tau,F)$ satisfies condition \req{nesc-exist}.}\1
\indent Assertion \req{exist6} follows from \rth{precisetr} (i) and \rdef{precisetr}.
Since $u\ind{reg}$ is \gsmod,
it follows, by \rth{gsmod-tr}, that $\tau\in \BBM_q(\bdw)$.

By \rth{reg-tr} (vi), $u\app{\CR(u)}u\ind{reg}$. Therefore  $w:=u\ominus u\ind{reg}$ vanishes on
$\CR(u)$ so that $w\leq U_F$. Note that $u- u\ind{reg}<w$ and therefore
\begin{equation}\label{exist8}
 u\leq u\ind{reg}\oplus w\leq U.
\end{equation}
By their definitions $\CS_0(u)=F_\tau$ and by \rth{precisetr} (vi) and \rlemma{bqF},
\begin{equation}\label{exist9}\BAL
  \CS(U)&=\CS(u\ind{reg})\cup \CS(U_F)= \CS(u\ind{reg})\cup b_q(F),\\
&=\CS_0(u)\cup b_q(F)=F_\tau\cup b_q(F).
\EAL\end{equation}

On the other hand, $\CR(U)\supset \CR(u\ind{\CR})=\CR(u)$ and, as $u\leq U$,
$\CR(U)\sbs \CR(u)$. Hence $\CR(U)=\CR(u)$ and $\CS(U)=\CS(u)$. Therefore, by \req{exist6} and
\req{exist9},
$F=\CS(U)= F_\tau\cup b_q(F)$. Thus $(\tau,F)$ satisfies \req{nesc-exist} and $\tr^c U=(\tau,F)$.
The fact that $U$ is the maximal solution with this trace follows from \req{exist8}.

The solution $U$ is \gsmod because both $u\ind{\CR}$ and $U_F$ are \gsmod solutions.
This fact, \wrto $U_F$, was proved in \cite{MV4}.

The uniqueness of the solution in the class of \gsmod solutions follows from \rth{gsmod-tr} (iv).

Finally we prove:\1
\num{\bf D}{\em If $\nu\in\BBM_q(\bdw)$ then the couple $(\tau,F)$ defined by
\begin{equation}\label{exist10}
 v:=\sup\set{u_{\nu\chr{E}}: E\in \CE_\nu},\;\tau:=\tr v,\; F:=\bdw\sms \CR(v)
\end{equation}
satisfies \req{nesc-exist}. This is the unique couple in $\GTC(\bdw)$ satisfying $\nu=\BBT(\tau,F)$.}\2
\indent The solution $v$ is \gsmod so that $\tau\in \BBM_q(\bdw)$.

We claim that $u:=v\oplus U_F$ is a solution with boundary trace $\tr^c u=(\tau,F)$.
Indeed $u\geq v$ so that $\CR(u)\sbs \CR(v)$.
On the other hand, since $\tau$ is $q$-locally finite on $\CR(v)=\bdw\sms F$, it follows
that $\CS(u)\sbs F$. Thus $\CR(v)\sbs \CR(u)$ and we conclude that $\CR(v)= \CR(u)$ and
$F=\CS(u)$. This also implies that $v=u\ind{reg}$.

Finally
$$\CS(u)=\CS(v)\cup\CS(U_F)=F_\tau\cup b_q(F),$$
so that $F$ satisfies \req{nesc-exist}.

The fact that, for $\nu\in\BBM_q(\bdw)$, the couple $(\tau,F)$ defined by \req{exist10} is the
only one in $\GTC(\bdw)$
satisfying $\nu=\BBT(\tau,F)$ follows immediately from the definition of these spaces.

Statements {\bf A--D} imply  (i)--(iv).
\eproof
\Remark  If $\nu\in \BBM_q(\bdw)$ then  $G$ and  $v$ as defined in \req{existence1} have
the following alternative representation:
\begin{equation}\label{existenceA}
  v:=\sup\set{u_{\nu \chr{Q}}: \, Q\in \CF_\nu},\q  G=\bigcup_{\CF_\nu}Q=\bigcup_{\CE_\nu}E,
\end{equation}
\begin{equation}\label{existence B}
 \CF_\nu:=\set{Q : \, Q\textrm{ \qop }, \;\nu(Q)<\infty}.\2
\end{equation}
\indent To verify this remark we first observe that \rlemma{set-approx}
implies that if  $A$ is a \qop set then there exists an increasing \seq of $q$-quasi closed sets
$\set{E_n}$ \sth $A=\cup_1^\infty E_n$. In fact, in the notation of \req{q-open-app}, we may choose
$E_n=F_n\sms L$ where $L=A'\sms A$ is a set of capacity zero.

Therefore $$\bigcup_{\CD_\nu}D\sbs \bigcup_{\CF_\nu}Q \sbs \bigcup_{\CE_\nu}E=:H.$$
On the other hand, if $E\in \CE_\nu$ then  $\muCR(u)(\tl E)=\muCR(u)(E)=\nu(E)<\infty$ and, by
 \rth{reg-tr} (ix), $\tl E$ is regular, i.e.,
there exists a \qop regular set $Q$ \sth $E\sbsq Q$. Thus $H=\bigcup_{\CD_\nu}D$.

If $D$ is a \qop regular set then $D=\cup_1^\infty E_n$, where $\set{E_n}$
is an increasing \seq of $q$-quasi closed sets. \Consy,
$$u_{\nu \chr{D}}=\lim u_{\nu \chr{E_n}}.$$
Therefore
$$\sup\set{u_{\nu \chr{Q}}: \, Q\in \CD_\nu}\leq \sup\set{u_{\nu \chr{Q}}: \, Q\in \CF_\nu}
\leq \sup\set{u_{\nu \chr{E}}: \, E\in \CE_\nu}.$$ On the other
hand, if $E\in \CE_\nu$ then there exists a \qop regular set $Q$
\sth $E\sbsq Q$. \Consy we have equality.
\par Note that, in view of this remark, \rth{exist} is an immediate
\cons of \rth{existence}.\3
\noindent{\bf Acknowledgment.} Both authors were partially sponsored by an EC grant through
the RTN Program "Front-Singularities", HPRN-CT-2002-00274 and by the French-Israeli cooperation program through grant
No. 3-1352. The first author (MM) also wishes to acknowledge the support of the Israeli Science Foundation through grant
No. 145-05.

\end{document}